\documentclass{article}
\usepackage[utf8]{inputenc}
\usepackage[a4paper, total={6in, 8in}]{geometry}

\usepackage{hyperref,amssymb,amsmath,amsthm,enumerate,enumitem,mathbbol,subcaption}

\usepackage{xcolor}
\usepackage{comment}
\usepackage[capitalise]{cleveref}
\usepackage{pgfplots}
\usepackage{mdframed}
\usepackage{placeins}
\usepackage{verbatim}

\usetikzlibrary{positioning,arrows}
\pgfplotsset{compat=1.17}

\newtheorem{thm}{Theorem}[section]
\newtheorem{lemma}[thm]{Lemma}
\newtheorem{prop}[thm]{Proposition}

\newtheorem{defn}[thm]{Definition}

\newtheorem{remark}[thm]{Remark}

\renewcommand{\a}{\alpha}
\renewcommand{\b}{\beta}

\newcommand{\de}{\partial}
\newcommand{\e}{\epsilon}
\newcommand{\g}{\mathcal{\mu}}
\newcommand{\m}{\mathfrak{m}}
\newcommand{\rk}{\textnormal{rk}\, }
\newcommand{\se}{\subseteq}

\renewcommand{\v}{\mathrm{v}}

\newcommand{\A}{\mathbb{A}}

\newcommand{\D}{\mathcal{D}}
\newcommand{\F}{\mathbb{F}}
\renewcommand{\H}{\mathcal{H}}
\renewcommand{\L}{\mathcal{L}}
\newcommand{\Mcal}{\mathcal{M}}
\newcommand{\N}{\mathbb{N}}
\newcommand{\Nil}{\mathcal{N}}
\renewcommand{\O}{\mathcal{O}}
\renewcommand{\P}{\mathbb{P}}

\newcommand{\Qcal}{\mathcal{Q}}
\newcommand{\R}{\mathcal{R}}

\newcommand{\Z}{\mathbb{Z}}

\newcommand{\<}{\langle}
\renewcommand{\>}{\rangle}

\title{Elliptic Loops}
\author{Massimiliano Sala$^1$ \and Daniele Taufer$^2$}
\date{ \small $^1$University of Trento - \href{mailto: maxsalacodes@gmail.com}{maxsalacodes@gmail.com}\\%
    $^2$KU Leuven - \href{mailto: daniele.taufer@gmail.com}{daniele.taufer@gmail.com} \footnote{Corresponding author. Postal address: Celestijnenlaan 200A, B-3001 Leuven (Belgium)}\\[2ex]%
    \normalsize March 2023}

\begin{document}
 
\maketitle

\begin{abstract}
    Given a local ring $(R,\m)$ and an elliptic curve $E(R/\m)$, we define its \emph{elliptic loop} as the points of $\P^2(R)$ projecting to $E$ under the canonical modulo-$\m$ reduction, endowed with an operation that extends the curve's addition.
    While its subset of points satisfying the curve's Weierstrass equation is a group, this larger object is proved to be a power-associative abelian algebraic loop, which is seldom completely associative.
    When an elliptic loop has no points of order $3$, its affine part is obtained as a stratification of a one-parameter family of elliptic curves defined over $R$, which we call \emph{layers}.
    Stronger associativity properties are established when $\m^e$ vanishes for small values of $e \in \Z$.
    When the underlying ring is $R = \Z/p^e\Z$, the infinity part of an elliptic loop is generated by two elements, the group structure of layers may be established and the points with the same projection and same order possess a geometric description.
\end{abstract}




\section{Introduction}

 Elliptic curves have been proving to be a field of intense research and fruitful applications.
 They were born as smooth plane projective cubics with (at least) a rational point defined over a field, but they have been generalized over arbitrary base schemes as proper smooth curves with geometrically connected genus-one fibers, with a prescribed zero section.

 Their main interest resides in the possibility of endowing these objects with an abelian group structure \cite[Theorem 2.1.2]{ECOverSchemes}, so that they constitute the abelian varieties of minimal positive dimension.
 When the underlying structure is a ring with prescribed properties, these objects may be embedded in a projective plane and their addition law may be explicitly determined \cite{Lenstra}.
 Among the several research lines offered by these objects, their lifting has played a crucial role both from theoretical \cite{Norman} and applied \cite{SatohAraki,Smart} purview.

 In the literature, lifts of elliptic curves over fields are usually investigated under several strong (classical) assumptions, including: lifting to integral domains, having a prescribed characteristic and enjoying certain levels of topological separation (usually, Hausdorff).
 With such requirements, the resulting lifts may still be proven to be abelian varieties \cite{Norman}. 
 Moreover, there are typically multiple lifts of a given curve, which are considered individually.
 Apart from the canonical lift \cite{Deuring}, the others are often indistinguishable and their interplay has not been deeply investigated.
 
 In this work, we relax these assumptions by investigating large objects made of all the projective points reducing to a given base curve.
 They are globally non-associative, yet they contain the usual lifts as associative sections.
 This way, we are not selecting special lifts of the base curve, but we are considering all of them at once.
 
 Precisely, we study the projective points over a local ring $(R,\m)$ with $6 \in R^*$ that project to an elliptic curve of odd order defined by a short Weierstrass model over $R/\m$, namely we define the set
 \[
    \L_{A,B}(R) = \{ (X:Y:Z) \in \P^2(R) \ | \ X^3 + AXZ^2 + BZ^3 - Y^2Z \in \m \}.
 \]
 We call these objects \emph{elliptic loops}, as we prove they are power-associative abelian algebraic loops (Proposition \ref{prop:loop} and Proposition \ref{prop:PowAss}).
 This result follows by establishing a sufficient criterion to ensure triples associativity in terms of the rank of a prescribed matrix (Theorem \ref{thm:ass}).
 
 Elliptic loops are not the first occurrence of quasigroups arising from operations between the points of cubics \cite{Ethe}.
 Elliptic curves may also be obtained as sections of cubic hypersurfaces, which were proved to be commutative Moufang loops \cite{Manin} that are not necessarily groups \cite{Kanevsky}.
 The loops defined in the present paper have weaker associativity properties, since they do not need to satisfy the Moufang identities.
 However, they are usually larger objects that appear to be natural candidates for investigating the lifts of the given curve, as they are proved to have a notably ordered geometric structure.
 
 The affine part of elliptic loops is completely described in terms of a new family of elliptic curves, which are defined over $R$ by non-invertible linear combinations of the Weierstrass polynomial and its Hessian polynomial, namely for every $t \in \m$ by the equation
 \[
    x^3 + Axz^2 + Bz^3 - y^2z - t ( 3Ax^2z + 3xy^2 + 9Bxz^2 - A^2z^3 ) = 0.
 \]
 Those projective cubics are again abelian varieties, which we call \emph{layers} since in absence of points of order $3$ they are proved to stratify the affine part of $\L_{A,B}(R)$ (Proposition \ref{prop:AffineIntersection}).
 For $t=0$, this construction agrees with the standard definition of elliptic curves over the ring $R$.
 
 The infinity part of $\L_{A,B}(R)$ exhibits a different behavior: layers may have a non-trivial intersection at infinity when $R$ has zero-divisors (Proposition \ref{prop:InfIntersection}).
 Moreover, layers do not necessarily cover the whole infinity part: if the base ring is Hausdorff with a discrete valuation $\v$, we show that layers never contain non-zero points at infinity $(X:1:Z)$ with $\v(Z) \leq \v(X)$ (Lemma \ref{lem:infinityValuation}).

 If the base ring is Hausdorff and $\m$ has a non-zero integer generator (e.g. a Cohen ring), the loop at infinity behaves similarly to the additive group of $\m$ (Proposition \ref{prop:tech}).
 When $R$ is a discrete valued ring, this is consistent with the known structure of formal groups constructed over large powers of $\m$, which can be made explicit via the discrete logarithm map \cite[Theorem IV.6.4]{Silverman}.
 
 When $\m^{\Nil} \neq \< 0 \>$ for small values of $\Nil \in \Z_{>0}$, we can explicitly exhibit triples of non-associative points of elliptic loops (Section \ref{sec:MoreAssNO}).
 Conversely, if $\m^{\Nil} = \< 0 \>$ for small values of $\Nil \in \Z_{>0}$, then special families of triples may be proved associative (Section \ref{sec:lownil}).
 
 When $R$ is Hausdorff and $\m^{2} = \< 0 \>$, the points with any given finite torsion lying over the same base point $P$ may be obtained via the additive action of a group at infinity (Theorem \ref{thm:qTors}), therefore they are proved to be projectively collinear, and if $\m$ is a principal $\Z$-module, then they precisely are all the points of such line projecting to $P$ (Proposition \ref{prop:qLine}).

 When the underlying ring is $\Z/p^e\Z$, the infinity part of elliptic loops is generated (as a loop) by the points $(p:1:0)$ and $(0:1:p)$, which have both order $p^{e-1}$ (Theorem \ref{thm:InfStructure}). 
 These generators play different roles with respect to elliptic curves defined over the same ring, as the group $\langle (0:1:p) \rangle$ is proved to be a forbidden locus for layers (Lemma \ref{lem:ForbiddenInf}).
 Moreover, there are several other groups at infinity of the same order, namely the infinity part of layers, which is also proved to be a cyclic group of order $p^{e-1}$ (Proposition \ref{prop:LayStructure}).
 As a consequence, if we let $E$ be the projected elliptic curve over $\F_p$, every layer is isomorphic to either $\Z/p^{e-1}\Z \times E$ or to $\Z/p^{e}\Z$.
 
 \subsection{Paper organization}
 This paper is organized as follows.
 \begin{itemize} \itemsep-0.2em
     \item In Section \ref{sec:Preliminaries}, we recall the known definitions and results that we employ throughout the work.
     \item Section \ref{sec:ELoops} is devoted to defining elliptic loops and proving that they are abelian algebraic loops.
     \item Novel properties of the Hessian of a Weierstrass polynomial are presented in Section \ref{sec:Hessian}.
     \item The power-associativity of elliptic loops is proved in Section \ref{sec:powass}.
     \item In Section \ref{sec:Layers}, we define layers and we prove that they stratify the affine points of elliptic loops.
     \item Several properties of the infinity loop are examined in Section \ref{sec:InfinityPart}.
     From this section, the base ring is often required to be Hausdorff.
     \item In Section \ref{sec:MoreAss}, we exhibit further associativity properties: although elliptic loops are proved to be almost never associative, weak forms of associativity are established under certain conditions.
     \item When the above conditions are met, we characterize in Section \ref{sec:PointsFiniteOrd} the geometry of points with order dividing a prescribed integer.
     \item Noteworthy structure results hold when the underlying ring is $R=\Z/N\Z$, which is examined in Section \ref{sec:ZNZ}.
     \item Finally, conclusions and further work are proposed in Section \ref{sec:conclusions}.
 \end{itemize}

 The Magma \cite{Magma} code for verifying the symbolic computations involved in this work may be found in \cite{MyCode}.
 For the sake of readability, the several formal verifications involved in the proofs of the current paper are omitted, but the reader can straightforwardly verify them by running the corresponding Magma code, which contains precise references to the related results.

 
 \section{Preliminaries} \label{sec:Preliminaries}
 
 In this paper, $R$ is always a local commutative ring with unity and $6 \in R^*$, and its maximal ideal is denoted by $\m$.
 We let $\pi:R \to R/\m$ be the canonical projection into the residue field, and to simplify the notation we also denote every componentwise or coefficientwise projection by $\pi$, i.e. for every $n \in \Z_{\geq 1}$ we write
 \[
 \pi\big( (X_i)_{1 \leq i \leq n} \big) = \big( \pi(X_i) \big)_{1 \leq i \leq n} \in (R/\m)^n \quad {\rm and } \quad  \pi\left( \sum_{d = 0}^n f_d x^d\right) = \sum_{d = 0}^n \pi(f_d) x^d \in (R/\m)[x].
 \]
 

 Given some elements $u_1,\dots,u_s$ of a ring $\R$, we denote by $\< u_1,\dots,u_s \>$ the ideal they generate in $\R$.
 In the context of the present paper, $\R$ will be either the local ring $R$ or a polynomial ring with coefficients in $R$.
 
 A tuple in $(X_0, \dots, X_n) \in R^{n+1}$ is called \emph{primitive} if $\<X_0, \dots, X_n\> = R$, which for local rings is equivalent to having at least one entry $X_i$ not in $\m$.
 The group of units $R^* = R \setminus \m$ acts on the primitive $(n+1)$-tuples by componentwise multiplication, and its quotient by this action is the \emph{projective $n$-space} over $R$, denoted by $\P^n(R)$.
 The orbit of an element $(X_0, \dots, X_{n})$ is denoted by $(X_0 : \dots : X_{n})$.
 We say that a point is \emph{affine} if its last coordinate does not belong to $\m$, otherwise we say it is \emph{at infinity}.
 
 Let $n,m \in \Z_{\geq 1}$ and $A \in M_{n \times m}(R)$.
 For every integer $1 \leq t \leq \min\{n,m\}$, the \emph{t-minor ideal} $I_t(A)$ is defined as the ideal generated by the $t$-th minor determinants of $A$, namely the ideal generated by all the determinants of the $t \times t$ submatrices of $A$.
 By convention we set $I_0(A) = R$, and for every $t > \min\{n,m\}$ we set $I_t(A) = \<0\>$.
 We also recall \cite[Chapter 4]{Brown} that the rank of $A$ is defined by
 \[
    \rk A = \max_{t \in \Z_{\geq 0}} \{ I_t(A) \neq \<0\> \}.
 \]
 This notion of rank over rings provides us with a practical way of testing if two projective points are equal \cite[Lemma 6]{SalTau}:
 \[
    (X_0:\dots:X_n) = (Y_0:\dots:Y_n) \in \P^n(R) \quad \iff \quad \rk \begin{bmatrix}
            X_0 & \dots & X_n \\
            Y_0 & \dots & Y_n
    \end{bmatrix} = 1.
 \]
 Given a homogeneous polynomial $F \in R[x_0,\dots,x_n]$, its evaluation on a given projective point $P = (X_0:\dots:X_n) \in \P^{n}(R)$ is unique modulo $R^*$, and we equivalently denote it by $F(P)$ or $F(X_0,\dots,X_n)$.
 
 In this work, we only consider the \emph{projective plane} $\P^2(R)$.
 Irreducible projective plane cubics with a rational point are called \emph{elliptic curves}.
 Since $6 \in R^*$, it is well-known \cite[Section III.1]{Silverman} that, up to a change of coordinates, they all arise from zero-sets of short Weierstrass polynomials, i.e. there are $A,B \in R$ such that 
 \[
    \Delta_{A,B} = -(4A^3+27B^2) \in R^*,
 \]
 and such that the elliptic curve is given by
 \[
    E_{A,B}(R) = \{ (X:Y:Z) \in \P^2(R) \ | \ X^3+AXZ^2+BZ^3-Y^2Z = 0 \}.
 \]
 
 Lenstra has proved \cite[Section 3]{Lenstra} that these objects have a group structure, which arises from a combination of different addition laws \cite[Section 4]{BosLen}.
 Whenever the projected curve $E_{A,B}(R/\m)$ has no $2$-torsion points, this operation may be entirely described by the addition law corresponding to $(0:1:0)$.
 Since this is a slight assumption in the computational practice (e.g. for ordinary curves employed in cryptographic applications), we assume that it is always the case.
 We recall here for convenience the explicit formulation of this addition law, with the minor correction pointed out in \cite[Section 2.11]{Washington}:
 \begin{equation*}
     (X_1:Y_1:Z_1) +_{(0:1:0)} (X_2:Y_2:Z_2) = (T_1:T_2:T_3),
 \end{equation*}
 where
 \begin{align*}
    T_1    = \ & Y_1Y_2(X_1Y_2 + X_2Y_1) - AX_1X_2(Y_1Z_2 + Y_2Z_1) - A(X_1Y_2 + X_2Y_1)(X_1Z_2 + X_2Z_1) \\
    & - 3B(X_1Y_2 + X_2Y_1)Z_1Z_2 - 3B(X_1Z_2 + X_2Z_1)(Y_1Z_2 + Y_2Z_1) + A^2(Y_1Z_2 + Y_2Z_1)Z_1Z_2, \\[0.1cm]
    T_2    = \ & Y_1^2Y_2^2 + 3AX_1^2X_2^2 + 9BX_1X_2(X_1Z_2 + X_2Z_1) - A^2X_1Z_2(X_1Z_2 + 2X_2Z_1) \\
    & - A^2X_2Z_1(2X_1Z_2 + X_2Z_1) - 3ABZ_1Z_2(X_1Z_2 + X_2Z_1) - (A^3 + 9B^2)Z_1^2Z_2^2, \\[0.1cm]
    T_3    = \ & 3X_1X_2(X_1Y_2 + X_2Y_1) + Y_1Y_2(Y_1Z_2 + Y_2Z_1) + A(X_1Y_2 + X_2Y_1)Z_1Z_2 \\
    & + A(X_1Z_2 + X_2Z_1)(Y_1Z_2 + Y_2Z_1) + 3B(Y_1Z_2 + Y_2Z_1)Z_1Z_2.
\end{align*}

The following lemma provides us with a clean and computationally efficient way of evaluating the above addition law.

\begin{lemma} \label{lem:AddLaw}
    Let $Q_1,Q_2,Q_3,Q_4 \in R[X_i,Y_i,Z_i]_{i\in\{1,2\}}$ be the homogeneous polynomials of bidegree $(1,1)$ defined by
    \begin{align*}
        Q_1 &= -A X_1 Z_2 - A X_2 Z_1 - 3 B Z_1 Z_2 + Y_1 Y_2,
        && Q_2 = A^2 Z_1 Z_2 - A X_1 X_2 - 3 B X_1 Z_2 - 3 B X_2 Z_1, \\
        Q_3 &= A Z_1 Z_2 + 3 X_1 X_2,
        && Q_4 = A X_1 Z_2 + A X_2 Z_1 + 3 B Z_1 Z_2 + Y_1 Y_2.
    \end{align*}
    By defining
    \begin{align*}
        T_1 &= (X_1 Y_2 + X_2 Y_1) Q_1 + (Z_1 Y_2 + Z_2 Y_1) Q_2, \\
        T_2 &= Q_1 Q_4 - Q_2 Q_3, \\
        T_3 &= (X_1 Y_2 + X_2 Y_1) Q_3 + (Z_1 Y_2 + Z_2 Y_1) Q_4,
    \end{align*}
    we have
    \[
        (T_1:T_2:T_3) = (X_1:Y_1:Z_1) +_{(0:1:0)} (X_2:Y_2:Z_2).
    \]
\end{lemma}
\proof Straightforward computation.
\endproof

 \subsection{Discrete valued ring}
 
 Let $\Z \cup \{\infty\}$ be the ordered additive monoid obtained from $(\Z,+,\leq)$ by adding the element $\infty$, such that for every $\gamma \in \Z \cup \{\infty\}$ it satisfies
 \[
    \infty + \gamma = \infty, \quad \text{and} \quad \infty \geq \gamma.
 \]
 When the map
 \[
    \v : R \to \Z \cup \{\infty\}, \quad X \mapsto \v(X) = \sup_{i \in \Z} \{ X \in \m^i \},
 \]
 satisfies, for every $X,Y \in R$,
 \[
    \v(XY) = \v(X)+\v(Y), 
 \]
 it is called the \emph{$\m$-adic valuation} of $R$, and it turns this ring into a \emph{discrete valued ring}, which corresponds to the classical \emph{discrete valuation ring} (DVR) when $R$ is a domain.
 Sometimes $\v$ is required to be surjective \cite{Manis}, but we remark that this is not the case here, as we are mainly interested in proper subsets of $\Z \cup \{\infty\}$.
 
 It is easy to see that $R^*=\{X \in R \ | \ \v(X) = 0\}$, and $\m = \{X \in R \ | \ \v(X) \geq 1 \}$.
 Moreover, for every $i \in \Z_{\geq 0}$ the ideal $\m^i$ is $\v$-closed, namely 
 \[
    X \in \m^i, \ \v(Y) \geq \v(X) \implies Y \in \m^i.
 \]
 Finally, for any integer $n \in \Z$, we denote the valuation of $n$-times the unity of $R$ by $\v(n) = \v(n\cdot 1_R)$.

 \subsection{Hausdorff rings}
 
 Many results about points at infinity will require the considered ring $R$ to be \emph{Hausdorff}, namely its $\m$-topology is Hausdorff, i.e.
 \[
    \bigcap_{i \in \Z_{\geq 0}} \m^i = \<0\>.
 \]
 For Hausdorff rings we may define the \emph{nilpotency} of $\m$ as
 \[
    \Nil = \inf_{i \in \Z_{\geq 1}} \big\{ \m^i = \<0\> \big\}.
 \]
 We remark that a Hausdorff ring may have $\m^i \neq \<0\>$ for every $i \in \N$, e.g. if $R = \Z_p$ is the ring of $p$-adic integers. In such cases, we have $\Nil = \inf \emptyset = \infty$.
 
 In Hausdorff rings, the inclusion sequence
 \[
    R = \m^{0} \supset \m \supset \dots \supset \m^{\Nil} = \<0\>
 \]
 is strictly decreasing until it stabilizes to $\<0\>$.
 
 If $R$ is a Hausdorff discrete valued ring, then $0$ is the unique element of maximal (infinite) valuation.
 There are many examples of such rings.
 As an instance, when the maximal ideal is principal $\m = \<\g\>$, we call its generator $\g$ the \emph{uniformizer} of $R$, and it is well known that the following are equivalent.
\begin{itemize}
    \item $R$ is Noetherian,
    \item $R$ is Hausdorff,
    \item either $R$ is a DVR, or there exists $n \in \Z_{\geq 1}$ such that the ideals of $R$ are precisely $\{\m^i\}_{0 \leq i \leq n}$.
\end{itemize}

The ring $R$ is a DVR precisely when $\Nil = \infty$.
If $\Nil$ is finite, it coincides with the $n$ of the above equivalence.
Moreover, $R$ is a field precisely when $\Nil = 1$, and in this case the $\m$-adic valuation of $R$ is the trivial valuation, namely for every $X \neq 0$ we have $\v(X) = 0$.
When this is not the case, for every non-zero element $X \in R$ there are uniquely determined $u \in R^*$ and $e_X \in \{0,\dots,\Nil-1\}$ such that
\[
    X = u \g^{e_X}.
\]
In this case, the $\m$-adic valuation of $R$ is the map $\v(X) = e_X$.

 \subsection{Finite rings}
 
 When the ring $R$ is local and finite, we have convenient enumerative properties \cite{GaMcDO}, which we recall here: 
 $\m$ is a nilpotent ideal, the residue field is a finite field $R/\m \simeq \F_q$ for some prime power $q = p^t$, and every $\m^{i}/\m^{i+1}$ is a finite-dimensional $\F_q$-vector space, so there is a positive $e \in \Z_{>0}$ such that
 \[
    |R| = q^e, \quad |\m| = q^{e-1}.
 \]
 
 One can explicitly count the points in projective spaces by observing that affine elements may be assumed to have the last entry equal to $1$, while the infinity part of $\P^n(R)$ may be realized by $|\m|$ copies of $\P^{n-1}(R)$. Hence, we inductively have
 \[
    |\P^n(R)| = \sum_{i = 0}^n |R|^{n-i}|\m|^{i} = q^{en} \left( 1 + \frac{1}{q} + \dots + \frac{1}{q^n} \right) = \frac{q^{en+1}-(q^{n})^{e-1}}{q-1}.
 \]

 \subsection{Non-associative structures}
 We recall some basic definitions and results about non-associative algebraic structures \cite{Schafer}.\\
 A \emph{magma} $\Mcal$ is a set equipped with a (closed) binary operation $+: \Mcal \times\Mcal \to \Mcal$.\\
 A \emph{quasigroup} $\Qcal$ is a magma satisfying the \emph{Latin square property}, i.e.
 \[ \forall P,Q \in \Qcal \ \exists !\ X, Y \in \Qcal \ : \ \begin{cases}
     P + X = Q, \\
     Y + P = Q.
 \end{cases} \]
 A \emph{loop} $\L$ is a quasigroup with \emph{identity}, i.e. there is $\O \in \L$ such that
 \[ \forall P \in \L \ : \ P + \O = \O + P = P. \]
 A magma (resp. quasigroup, loop) is called \emph{abelian} if its operation $+$ is commutative.
 
 An abelian loop $\L$ is called
 \begin{itemize}
     \item \emph{alternative loop}, if $\forall P,Q \in \L \ : \ P + (P + Q) = (P + P) + Q$,
     \item \emph{Jordan loop}, if $\forall P,Q \in \L \ : \ (P + P) + (P + Q) = P + \big(Q + (P + P)\big)$,
     \item \emph{Moufang loop} (or \emph{Bol loop}), if $\forall P,Q,R \in \L \ : \ \big(P + (Q + R)\big) + R = \big((P + R) + R\big) + Q$,
     \item \emph{power-associative loop}, if the subloop generated by any element is associative,
     \item \emph{diassociative loop}, if the subloop generated by any subset of size at most two is associative.
 \end{itemize}
 Finally, we observe that diassociative clearly implies power-associative, while it can be shown that Moufang implies both alternative and dissociative.
 
 
 \section{Elliptic Loops} \label{sec:ELoops}
 
 \begin{defn}[Elliptic Loop]
  Let $A,B \in R$ be elements defining an elliptic curve $E_{A,B}(R/\m)$ without points of even order.
  We define its \emph{Elliptic Loop} as the set
  \[
    \L_{A,B}(R) = \{ (X:Y:Z) \in \P^2(R) \ | \ X^3 + AXZ^2 + BZ^3 - Y^2Z \in \m \},
  \]
  endowed with the operation $+_{(0:1:0)}$, simply referred to as $+$.
  When the coefficients and the underlying ring are understood, we simply denote it by $\L$.
  The polynomial $x^3+Axz^2+Bz^3-y^2z \in R[x,y,z]$ is referred to as the \emph{defining polynomial of} $\L$.
 \end{defn}

 \begin{remark}
     The operation $+$ always depends on the considered $A,B \in R$.
     Moreover, as recalled in the previous section, the absence of points of order $2$ in $E(R/\m)$ guarantees that the operation $+$ is always well-defined between points of $\L$. 
 \end{remark}
 
 Elliptic loops canonically project on their underlying elliptic curve via $\pi$, and it is easy to see that when $R$ is finite this is a $|\m|^2$-covering of the curve.
 We define their \emph{affine} (resp. \emph{infinity}) \emph{points} $\L^a$ (resp. $\L^{\infty}$) as those projecting to affine points of $E_{A,B}(\F_q)$ (resp. to $\O$), namely
 \[
    \L^{a} = \{ (X:Y:Z) \in \L \ | \ Z \not\in \m \}, \quad \L^{\infty} = \{ (X:1:Z) \}_{X,Z \in \m}.
 \]
 We remark that, as a set, $\L^{\infty} = \pi^{-1}(\O)$ is independent of the considered curve parameters $A,B \in R$, so it is common to all the elliptic loops over the same ring $R$, but its operation depends on the coefficients of the underlying curve.
 
 The elliptic loop constructed over $E = E_{A,B}(R/\m)$ contains (set-theoretically) the points of every possible lift of $E$ over $R$.
 The next lemma shows that its affine part is precisely made of the affine points of classical elliptic curves over $R$ that project on $E$.
 
 \begin{lemma} \label{lem:AffinePointsofLoop}
 Let $\L = \L_{A,B}(R)$ be an elliptic loop and $P \in \L^a$ be one of its affine points.
 For every $\a \in \m$ there exists $\b \in \m$ such that $P \in E_{A+\a,B+\b}(R)$.
 \end{lemma}
 \proof
 Since $P$ is affine, we may assume $P = (X:Y:1)$.
 For every $\a \in \m$, we define
 \[
    \b = Y^2 - X^3 - (A+\a)X - B.
 \]
 We immediately verify that $P$ satisfies the Weierstrass equation of $E_{A+\a,B+\b}(R)$. Since $\a X \in \m$ and by definition of elliptic loop $Y^2 - X^3 - AX - B \in \m$, then also $\b \in \m$.
 \endproof
 Unlike the affine component, the infinity part $\L^{\infty}$ usually contains more points than those arising from elliptic curves over $R$, such as $\{(0:1:\a)\}_{\a \in \m \setminus \{\O\}}$.

 We now show that the name is well-given, i.e. $\L$ is always an algebraic loop. 
 We begin by proving that it is a set with a well-defined commutative binary operation, namely an abelian magma.
 
 \begin{prop} \label{prop:magma}
 Every elliptic loop $\L = \L_{A,B}(R)$ is an abelian magma with identity $\O = (0:1:0)$.
 Moreover, every element $(X:Y:Z) \in \L$ has a unique inverse
 \[
    -(X:Y:Z) = (X:-Y:Z).
 \]
 \end{prop}
 \proof
 The addition law $+$ is symmetric in the two addenda, hence it is commutative.
 Moreover, it is defined by polynomial relations in their entries, hence it commutes with $\pi$, so for every $P_1, P_2 \in \L$ we have
 \[
    \pi(P_1 + P_2) = \pi(P_1) + \pi(P_2) \in E_{A,B}(R/\m).
 \]
 This shows that $P_1 + P_2 \in \L$, hence $+$ is a well-defined binary operation on $\L$.
 
 To show that $\O$ acts as identity, we explicitly compute
 \[
    (X:Y:Z) + (0:1:0) = (XY:Y^2:YZ).
 \]
 Since $E(R/\m)$ has no points of order $2$, the latter needs to be a proper projective point, hence we have $Y \not\in \m$ and we conclude $(XY:Y^2:YZ) = (X:Y:Z) \in \P^2(R)$. 
 
 To show that $(X:-Y:Z)$ is an inverse of $(X:Y:Z)$, we compute
 \[
    (X:Y:Z) + (X:-Y:Z) =
    (0 : 3AX^4 + 18BX^3Z - 6A^2X^2Z^2 - 6ABXZ^3 + Y^4 + (-A^3 - 9B^2)Z^4 : 0),
 \]
 and since the operation is well-defined, the result equals $\O$.
 
 As for the uniqueness part, let $P_1 = (X_1:Y_1:Z_1)$ and $P_2 = (X_2:Y_2:Z_2)$ be two points of $\L$ such that
 \[
    (T_1:T_2:T_3) = P_1 + P_2 = \O \in \P^2(R).
 \]
 We observe that, as polynomials in $R[X_i,Y_i,Z_i]_{i \in \{1,2\}}$, we have
 \[
    (X_1Y_2 + X_2Y_1) T_2 \in \< T_1,T_3 \>, \quad (Y_1Z_2 + Y_2Z_1) T_2 \in \< T_1,T_3 \>.
 \]
 hence $(T_1:T_2:T_3) = \O$ implies $X_1Y_2 + X_2Y_1 = Y_1Z_2 + Y_2Z_1 = 0$. Since the base curve has odd order, we have $Y_1,Y_2 \not\in \m$, hence we conclude
 \[
    (X_2:Y_2:Z_2) = \left( -\frac{Y_1}{Y_2}X_2 : -\frac{Y_1}{Y_2}Y_2 : -\frac{Y_1}{Y_2}Z_2 \right) = ( X_1 : -Y_1 : Z_1 ),
 \]
 i.e. the unique inverse of $(X_1:Y_1:Z_1)$ is $(X_1:-Y_1:Z_1)$.
 \endproof

 If elliptic loops were associative, they would be groups.
 Although this is almost never the case (Section \ref{sec:MoreAssNO}), a weak form of associativity always holds.
 
 \begin{lemma} \label{lem:WeakAssociativity}
 Let $P,Q \in \L$ be two points of an elliptic loop. Then
 \[
    P + (-P+Q) = Q.
 \]
 \end{lemma}
 \proof Let $Q = (Q_1:Q_2:Q_3)$. We can directly compute by means of an algebraic calculator
 \[
    (T_1:T_2:T_3) = P+(-P+Q),
 \]
 and we symbolically verify that
 \[
    I_2 \left( \begin{bmatrix}
     T_1 & T_2 & T_3 \\
     Q_1 & Q_2 & Q_3
    \end{bmatrix} \right) = \<0\>,
 \]
 which means $(T_1:T_2:T_3) = Q \in \P^2(R)$.
 \endproof
 
 From Lemma \ref{lem:WeakAssociativity} we can show that elliptic loops satisfy the Latin square property, so they are quasigroups with identity, i.e. algebraic loops.
 
 \begin{prop} \label{prop:loop}
 Every elliptic loop $\L$ is an abelian algebraic loop.
 \end{prop}
 \proof
 By Proposition \ref{prop:magma} we know that $\L$ is an abelian magma with identity, hence it is sufficient to show that for every $P,Q \in \L$ there is a unique $R \in \L$ such that $P+R=Q$.
 By Lemma \ref{lem:WeakAssociativity} a solution is $R = Q-P$, and it is unique because if $R_2 \in \L$ satisfies the same equation, by the same lemma we have
 \[
    R_2 = -P + ( P + R_2 ) = -P + Q,
 \]
 which implies $R=R_2$.
 \endproof

 One may check that even for relatively simple rings (such as $R = \Z/p^e\Z$, see Section~\ref{sec:ZNZ}) elliptic loops need not be alternative loops, Jordan loops or diassociative loops, hence neither Moufang loops.
 However, in Section \ref{sec:powass} we will prove that they are necessarily power-associative loops.
 
 
 \section{The Hessian of an elliptic curve} \label{sec:Hessian}
 
 Given a projective degree-$d$ hypersurface in $\P^n(R)$, we may construct its Hessian as the hypersurface of degree $(d-2)(n+1)$ defined by the vanishing of its Hessian matrix.
 We refer to \cite{CilOtt} for a general treatment of such surfaces.
 In this work we are primarily interested in the Hessian polynomial, i.e. the polynomial defining the Hessian hypersurface: given a homogeneous degree-$d$ polynomial $F \in R[x_0, \dots, x_n]$, its \emph{Hessian} (polynomial) is defined by
 \[
    \H_F = \det ( \de_{x_i}\de_{x_j}F )_{0 \leq i,j \leq n} \in R[x_0, \dots, x_n].
 \]

 It is easy to see that the surfaces defined by a polynomial and its Hessian lie in the same projective space precisely in the following cases:
 \begin{itemize}
     \item $d = 3, n = 2$: the curves defined by $F$ and $\H_F$ are plane projective cubics.
     \item $d = 4, n = 1$: they are quartics defined on a projective line.
 \end{itemize}
 
 Here we only consider the first case: the Hessian polynomial of a short Weierstrass polynomial $F = x^3+Axz^2+Bz^3-y^2z \in R[x,y,z]$ is
 \[
    \H_F = -8( 3Ax^2z + 3xy^2 + 9Bxz^2 - A^2z^3 ) \in R[x,y,z].
 \]
 As we are interested in the zeroes of $\H_F$ over rings with $2 \in R^*$, we can neglect the factor $-8$.
 We also recall that the Hessian of an elliptic curve of non-zero $j$-invariant defines another elliptic (hence non-singular) curve \cite[Proposition 5.11]{CilOtt}.
 
 Regardless of their smoothness, the curves defined by $F$ and $\H_F$ over $R$ lie in the same ambient space, so we may investigate their intersection properties.
 In particular, the following lemma extends \cite[Exercise III.3.9]{Silverman}, and it shows that they may intersect only on lifts of $3$-torsion points.
 
 \begin{lemma} \label{lem:3tors}
 Let $P \in \L$ be a point of the elliptic loop with defining polynomial $F$.
 Then
 \[
    3 \pi(P) = \O \quad \iff \quad \H_F(P) \in \m.
 \]
 \end{lemma}
 \proof
 We formally compute the polynomials $D_1,D_2,D_3,T_1,T_2,T_3 \in R[x,y,z]$ such that
 \[
    (D_1:D_2:D_3) = (x:y:z) + (x:y:z), \quad 
    (T_1:T_2:T_3) = (x:y:z) + (D_1:D_2:D_3).
 \]
 We straightforwardly verify that
 \[
    T_2 \H_F \in \< F, D_2 z + D_3 y, D_1 z - D_3 x, D_1 y + D_2 x \>.
 \]
 If $3 \pi(P) = \O$, then $T_2(P) \in R^*$ and $\pi(2P) = \pi(-P) \in \P^2(R/\m)$, hence for each $P = (X:Y:Z)$ we have
 \[
    I_2 \left(\begin{bmatrix}
        D_1(P) & D_2(P) & D_3(P) \\
        X & -Y & Z
    \end{bmatrix} \right) \se \m.
 \]
 Moreover, since $P \in \L$ then $F(P) \in \m$, therefore we conclude $\H_F(P) \in \m$.
 
 On the other side, we symbolically verify the ideal inclusion
 \[
    \< T_1,T_3 \> \se \< F,\H_F \>.
 \]
 Since $F(P) \in \m$ by definition of $\L$, then $\H_F(P) \in \m$ implies $T_1(P),T_3(P) \in \m$, from which we conclude that $3\pi(P) = \O$.
 \endproof
 
 We now observe that the addition law $+$ respects linear combinations of $F$ and $\H_F$.
 
 \begin{lemma} \label{lem:HessianCombinations}
 Let $F$ be the defining polynomial of an elliptic loop $\L$, and $\a,\b \in R$.
 For every $P_1,P_2 \in \P^2(R)$ such that
 \[
    ( \a F + \b \H_F )(P_1) = 0 = ( \a F + \b \H_F )(P_2),
 \]
 we have
 \[
    ( \a F + \b \H_F )(P_1 + P_2) = 0.
 \]
 \end{lemma}
 \proof We define $T_{\a,\b} = \a F + \b \H_F \in R[x,y,z]$ and we symbolically verify that in the polynomial ring $R[x_i,y_i,z_i]_{i \in \{1,2\}}$ we have
 \begin{equation} \label{eq:Tincl}
    T_{\a,\b}\big( (x_1:y_1:z_1) + (x_2:y_2:z_2) \big) \in \< T_{\a,\b}(x_1,y_1,z_1), T_{\a,\b}(x_2,y_2,z_2) \>.
 \end{equation}
 Thus, for every pair of points $P_1,P_2 \in \P^2(R)$ such that $T_{\a,\b}(P_1) = 0 = T_{\a,\b}(P_2)$, we always have $T_{\a,\b}(P_1 + P_2) = 0$.
 \endproof
 
 \begin{remark}
 Lemma \ref{lem:HessianCombinations} does not ensure that $P_1 + P_2$ is a proper projective point. As an instance, if the elliptic curve defined by $F$ has $2$-torsion points, this operation may produce points with all entries in $\m$.
 However, this lemma may be extended to other addition laws defined in \cite{BosLen}, which may also be proved to projectively agree with a few more assumptions \cite{TauThesis}.
 This way, this result holds also if the elliptic curve defined by $F$ has points of even order.
 \end{remark}
 
 
 \section{Power-Associativity} \label{sec:powass}

 Even for non-associative loops, one may define a notion of \emph{multiple} for a given point $P$ recursively: for every $n \in \Z$, we set
 \[
    0 P = \O, \quad (-n)P = n(-P), \quad (n+1) P = n P + P.
 \]
 We denote the set of multiples of a given point $P$ by
 \[
    \langle P \rangle = \{ nP\}_{n \in \Z}.
 \]
 If the loop operation is associative when restricted to the multiples of every given point, then we say that the loop is \emph{power-associative}.
 Equivalently, the multiples of an element in power-associative loops form a group, i.e. for every $n,m \in \Z$ we have
 \[
    (n+m)P = nP + mP.
 \]
 Those loops are computationally friendly, since fast-multiplication techniques may be applied to evaluate $nP$ in at most $O(\log n)$ point additions.
 
 In this section, we prove that elliptic loops are power-associative.
 To show that, we employ the following objects.
 
 \begin{defn}[Associativity matrix]
    Let $\L$ be an elliptic loop. 
    Given $n$ points $P_1, \dots, P_n \in \L$, we define their \emph{associativity matrix} as
    \[
        \A(P_1,\dots,P_n) = \begin{bmatrix}
            F(P_1) & \dots & F(P_n) \\
            \H_F(P_1) & \dots & \H_F(P_n)
        \end{bmatrix} \in M_{2 \times n}(R).
    \]
 \end{defn}
 
 The following theorem shows that a sub-maximal rank of the associativity matrix is a sufficient condition for establishing points associativity.
 
 \begin{thm} \label{thm:ass}
 Let $P_1, P_2, P_3 \in \L$ be points of an elliptic loop.
 \begin{enumerate}[label=$(\roman*)$]
     \item \label{i} We have
     \[
        \rk \A(P_1,P_2+P_3) \leq \rk \A(P_1,P_2,P_3).
     \]
     \item \label{ii} If $\rk \A(P_1,P_2,P_3) \leq 1$, then the triple is associative, i.e.
     $(P_1 + P_2) + P_3 = P_1 + (P_2 + P_3)$.
 \end{enumerate}
 \end{thm}
 \proof  Let $F$ be the defining polynomial of $\L$.
 
 \ref{i} It is sufficient to show that for $i \in \{1,2\}$ we have $I_i\big(\A(P_1,P_2+P_3)\big) \se I_i\big(\A(P_1,P_2,P_3)\big)$.
 
 Since the $1$-minor ideal of a matrix is the ideal generated by its entries, the case $i=1$ follows from the inclusion \eqref{eq:Tincl} of Lemma \ref{lem:HessianCombinations}, by considering $T_{1,0}$ and $T_{0,1}$.
 As for the $2$-minors, we formally compute in $R[x_i,y_i,z_i]_{i\in\{1,2,3\}}$ the sum
 \[
    T = (x_2:y_2:z_2) + (x_3:y_3:z_3),
 \]
 and we explicitly verify that
 \[
    F(x_1,y_1,z_1) \H_F(T) - F(T) \H_F(x_1,y_1,z_1) \in I_2\Big( \A\big( (x_1:y_1:z_1), (x_2:y_2:z_2), (x_3:y_3:z_3)\big) \Big).
 \]
 We observe that $\big(F(x_1,y_1,z_1) \H_F(T) - F(T)\H_F(x_1,y_1,z_1)\big)(P_1,P_2,P_3)$ is the unique generator of $I_2\big(\A(P_1,P_2+P_3)\big)$, from which the inclusion between the $2$-minor ideals follows.
 
 \smallskip
 \ref{ii} We explicitly compute $S_1,S_2,S_3,T_1,T_2,T_3 \in R[x_i,y_i,z_i]_{i\in\{1,2,3\}}$ such that
 \begin{align*}
    (S_1:S_2:S_3) &= \big( (x_1:y_1:z_1) + (x_2:y_2:z_2) \big) + (x_3:y_3:z_3), \\
    (T_1:T_2:T_3) &= (x_1:y_1:z_1) + \big( (x_2:y_2:z_2) + (x_3:y_3:z_3) \big). 
 \end{align*}
 We straightforwardly verify that
 \[
    I_2\left( \begin{bmatrix} S_1 & S_2 & S_3 \\
    T_1 & T_2 & T_3
    \end{bmatrix}\right) \se I_2\Big( \A\big( (x_1:y_1:z_1), (x_2:y_2:z_2), (x_3:y_3:z_3)\big) \Big),
 \]
 hence when $\rk \A(P_1,P_2,P_3) \leq 1$ we have $I_2\big(\A(P_1,P_2,P_3)\big) = \< 0 \>$, so we conclude 
 \begin{align*}
    (P_1+P_2)+P_3 &= \big(S_1(P_1,P_2,P_3):S_2(P_1,P_2,P_3):S_3(P_1,P_2,P_3)\big)\\
    &= \big(T_1(P_1,P_2,P_3):T_2(P_1,P_2,P_3):T_3(P_1,P_2,P_3)\big) =
    P_1+(P_2+P_3),
 \end{align*}
 namely the addition law on $P_1,P_2,P_3$ is associative.
 \endproof
 
 By means of Theorem \ref{thm:ass}, we show that every elliptic loop is power-associative.
 
 \begin{prop} \label{prop:PowAss}
 Let $P \in \L$ be a point of an elliptic loop. Then $\langle P \rangle$ is a group.
 \end{prop}
 \proof Let $F$ be the defining polynomial of $\L$. By part \ref{ii} of Theorem \ref{thm:ass} it is sufficient to show that for every $n,m \in \Z$ we have
 \begin{equation} \label{eq:rkineq}
     \rk \begin{bmatrix}
    F(mP) & F(nP) \\
    \H_F(mP) & \H_F(nP)
    \end{bmatrix} \leq 1.
 \end{equation}
 Since $F(\O) = \H_F(\O) = 0$, $F(P) = F(-P)$ and $\H_F(P) = \H_F(-P)$, we may assume without losing of generality that $n \geq m \geq 1$.
 Moreover, we observe that the case $m=n$ is trivial, since 
 \begin{equation} \label{eq:m=n}
     I_2\left( \begin{bmatrix} F(nP) & F(nP) \\
    \H_F(nP) & \H_F(nP)
    \end{bmatrix}\right) = \< F(nP)\H_F(nP) - \H_F(nP)F(nP) \> = \< 0 \>.
 \end{equation}
 
 We prove by extended induction on $n \in \Z_{\geq 1}$ that the inequality \eqref{eq:rkineq} holds for every $1 \leq m \leq n$. \\
 $[n = 1]$: There is only one possible $m$, that is $m = n = 1$. \\
 $[1, \dots, n-1 \to n]$: By part \ref{i} of Theorem \ref{thm:ass} we have
 \[
    \rk \begin{bmatrix}
    F(mP) & F(nP) \\
    \H_F(mP) & \H_F(nP)
    \end{bmatrix} \leq
    \rk \begin{bmatrix}
    F(mP) & F\big((n-1)P\big) & F(P) \\
    \H_F(mP) & \H_F\big((n-1)P\big) & \H_F(P)
    \end{bmatrix}.
 \]
 When $m=n$ the inequality \eqref{eq:rkineq} holds by \eqref{eq:m=n}, while if $m<n$ then by inductive hypothesis
 \[
    {\rm rk} \begin{bmatrix}
    F(mP) & F\big((n-1)P\big) \\
    \H_F(mP) & \H_F\big((n-1)P\big)
    \end{bmatrix} \leq 1, \
    {\rm rk} \begin{bmatrix}
    F\big((n-1)P\big) & F(P) \\
    \H_F\big((n-1)P\big) & \H_F(P)
    \end{bmatrix}\leq 1, \
    {\rm rk} \begin{bmatrix}
    F(mP) & F(P) \\
    \H_F(mP) & \H_F(P)
    \end{bmatrix}\leq 1.
 \]
 Thus, we obtain
 \[
 \rk \begin{bmatrix}
    F(mP) & F\big((n-1)P\big) & F(P) \\
    \H_F(mP) & \H_F\big((n-1)P\big) & \H_F(P)
    \end{bmatrix} \leq 1,
 \]
 which concludes the inductive step.
 \endproof

 \section{Layers} \label{sec:Layers}
 
 \begin{defn}
    Let $\L$ be an elliptic loop with defining polynomial $F$. For every $t \in \m$, we define the $t$-\emph{layer} of $\L$ as
    \[
        L_t = \{ P \in \P^2(R) \ | \ (F-t\H_F)(P) = 0  \}.
    \]
 \end{defn}
 
 By definition $t \in \m$, then $L_t \se \L$.
 Moreover, the elliptic curve $E_{A,B}(R)$ is the $0$-layer of $\L_{A,B}(R)$, hence there is a projective abelian variety lying inside such a loop.
 The following proposition shows that this actually holds for every layer.
 
 \begin{prop}
 Let $\L$ be an elliptic loop and $t \in \m$. Then $L_t$ is a subloop of $\L$, which is a group with the addition law of $\L$.
 \end{prop}
 \proof The addition law of $\L$ is closed on $L_t \se \L$ by Lemma \ref{lem:HessianCombinations}, and they have the same unity $\O \in L_t$, so we only need to show that this operation is associative on $L_t$.
 Let $P_1,P_2,P_3 \in L_t$, then for every $i \in \{1,2,3\}$ we have
 \[
    F(P_i) = t \H_F(P_i).
 \]
 Thus, we have
 \[
    \A(P_1,P_2,P_3) = \begin{bmatrix}
        t \H_F(P_1) & t \H_F(P_2) & t \H_F(P_3) \\
        \H_F(P_1) & \H_F(P_2) & \H_F(P_3) 
    \end{bmatrix},
 \]
 which has vanishing $2$-minors. Hence, $\rk \A(P_1,P_2,P_3) \leq 1$, so $P_1,P_2,P_3$ associate by part \ref{ii} of Theorem \ref{thm:ass}.
 \endproof
 
 We denote the affine (resp. infinity) points of $L_t$ by
 \[
    L_t^a = L_t \cap \L^a, \quad L_t^{\infty} = L_t \cap \L^{\infty}.
 \]
 
 As shown by the following proposition, in absence of $3$-torsion points, the affine part of an elliptic loop is stratified by the affine parts of its layers.
 
 \begin{prop} \label{prop:AffineIntersection}
    Let $\L$ be an elliptic loop without points of order $3$.
    Then $\L^a$ is given by the disjoint union
    \[
        \L^a = \bigsqcup_{t \in \m} L_t^a.
    \]
 \end{prop}
 \proof Since $\L$ does not have points of order $3$, by Lemma \ref{lem:3tors} for every $P \in \L^a$ we have $\H_F(P) \in R^*$.
 Therefore, we may define
 \[
    n_P = F(P) \big(\H_F(P)\big)^{-1} \in \m,
 \] 
 so we have $P \in L_{n_P}$.
 As for the disjointness, we observe that $P \in L_s \cap L_t$ implies
 \[
    (s-t)\H_F(P) = 0,
 \]
 but $\H_F(P) \in R^*$, so we conclude $s=t$.
 \endproof
 
 We observe that for Weierstrass polynomials defining ordinary curves of cryptographic interest, the absence of $3$-torsion points is a standard assumption.
 
 The group structure of layers $L_t \se \L$ depends on the exact sequence of groups
 \begin{equation} \label{eq:ses}
     0 \to L_t^{\infty} \to L_t \xrightarrow{\pi} E_{A,B}(R/\m) \to 0,
 \end{equation}
 which generalizes the usual reduction sequence of elliptic curves \cite[Proposition VII.2.1]{Silverman}.
 Determining the group structure of $L_t^{\infty}$ may be challenging, depending on the properties of the underlying ring.
 For instance, in Section \ref{sec:ZNZ} we prove that this group is cyclic if the underlying ring is $\Z/p^e\Z$, but this result heavily relies on the special properties of this ring.
 
 
\section{The infinity part} \label{sec:InfinityPart}

In this section, we outline some general properties of the infinity parts of elliptic loops.

\begin{prop} \label{prop:tech}
Let $\L$ be an elliptic loop and $(X_1:1:Z_1), (X_2:1:Z_2) \in \L^{\infty}$ be two of its points adding to
\[
    (X:1:Z) = (X_1:1:Z_1) + (X_2:1:Z_2).
\]
If $e \in \Z_{\geq 1}$ is such that $X_1,Z_1,X_2,Z_2 \in \m^e$, then the following hold.
\begin{enumerate}[label=$(\roman*)$]
    \item \label{I} We have
    \[
        X \equiv X_1 + X_2 \bmod \m^{3e}, \quad Z \equiv Z_1 + Z_2 \bmod \m^{3e}.
    \]
    \item \label{II} For every $f \geq e$ and $\delta_X,\delta_Z \in \m^f$, if
    \[
        (X_m:1:Z_m) = (X_1+\delta_X:1:Z_1+\delta_Z) + (X_2:1:Z_2),
    \]
    then
    \[
        X_m \equiv X+\delta_X \bmod \m^{f+2e}, \quad Z_m \equiv Z+\delta_Z \bmod \m^{f+2e}.
    \]
    
    \item \label{III} If $R$ is a discrete valued ring and $n \in \Z$ is an integer with $\v(n) \leq 2e$, then for every $f \in \Z_{\geq 1}$ the $n^f$-th multiple
    \[
        (X_f:1:Z_f) = n^f(X_1:1:Z_1),
    \]
    satisfies
    \[
        X_f \equiv n^fX_1 \bmod \m^{3e+(f-1)\v(n)}, \quad Z_f \equiv n^fZ_1 \bmod \m^{3e+(f-1)\v(n)}.
    \]
    
    \item \label{IV} If $R$ is Hausdorff and it has a non-zero integer uniformizer, then for every $\a \in \Z$ we have
    \[
        (X_{\a}:1:Z_{\a}) = \a(X_1:1:Z_1),
    \]
    with
    \[
        X_{\a} \equiv \a X_1 \bmod \m^{3e + \v(\a) - 1}, \quad Z_{\a} \equiv \a Z_1 \bmod \m^{3e + \v(\a) - 1}.
    \]
\end{enumerate}
\end{prop}
\proof By Lemma \ref{lem:AddLaw} there are homogeneous polynomials $H_1,H_2,H_3 \in R[x_i,z_i]_{i\in\{1,2\}}$ of degrees $\deg(H_1) = \deg(H_3) = 3$ and $\deg(H_2) = 4$ such that
\[
    \rk\begin{bmatrix}
    X&1&Z\\
    X_1+X_2+H_1(X_1,Z_1,X_2,Z_2) & 1+H_2(X_1,Z_1,X_2,Z_2) & Z_1+Z_2+H_3(X_1,Z_1,X_2,Z_2)
    \end{bmatrix} \leq 1,
\]
which implies
\begin{align*}
    \big(1+H_2(X_1,Z_1,X_2,Z_2)\big) X &= X_1+X_2+H_1(X_1,Z_1,X_2,Z_2),\\
    \big(1+H_2(X_1,Z_1,X_2,Z_2)\big)Z &= Z_1+Z_2+H_3(X_1,Z_1,X_2,Z_2).
\end{align*}

\ref{I} Since $X_1,Z_1,X_2,Z_2 \in \m^e$ and the $H_i$'s are homogeneous of degree at least $3$, the above equations imply that both $X - (X_1+X_2)$ and $Z - (Z_1+Z_2)$ belong to $\m^{3e}$.

\ref{II} Since the $H_i$'s are homogeneous of degree at least $3$ and $f \geq e$, the factors $\delta_X$ and $\delta_Z$ always appear multiplied by some element of $\m^{2e}$, hence for every $1 \leq i \leq 3$ we have
\begin{equation} \label{eq:Himod}
    H_i(X_1,Z_1,X_2,Z_2) \equiv H_i(X_1+\delta_X,Z_1+\delta_Z,X_2,Z_2) \bmod \m^{f+2e}.
\end{equation}
By using again Lemma \ref{lem:AddLaw} we have
\[
    \rk\begin{bmatrix}
    X_m&1\\
    X_1+\delta_X+X_2+H_1(X_1+\delta_X,Z_1+\delta_Z,X_2,Z_2) & 1+H_2(X_1+\delta_X,Z_1+\delta_Z,X_2,Z_2)
    \end{bmatrix} \leq 1
\]
which by \eqref{eq:Himod} implies
\[
    \big(1+H_2(X_1,Z_1,X_2,Z_2)\big) X_m \equiv X_1+\delta_X+X_2+H_1(X_1,Z_1,X_2,Z_2) \bmod \m^{f+2e}.
\]
Since $\delta_X \big(1 + H_2(X_1,Z_1,X_2,Z_2)\big) \equiv \delta_X \bmod \m^{f+2e}$, we conclude
\begin{align*}
    \big(1+H_2(X_1,Z_1,X_2,Z_2)\big) (X_m-\delta_X) &\equiv X_1+X_2+H_1(X_1,Z_1,X_2,Z_2) \\
    &\equiv \big(1+H_2(X_1,Z_1,X_2,Z_2)\big) X \bmod \m^{f+2e},
\end{align*}
from which $X_m \equiv X + \delta_X$ follows by inverting $1+H_2(X_1,Z_1,X_2,Z_2)$ modulo $\m^{f+2e}$.
The analogous result on $Z_m$ follows in the same way.

\ref{III} We do it by induction on $f \geq 1$.

[$f = 1$] It follows by applying $n$ times part \ref{I}.

[$f \to f+1$] Since $\L$ is power-associative (Proposition \ref{prop:PowAss}), we can compute $n^{f+1}(X_1:1:Z_1)$ as $n\big(n^{f}(X_1:1:Z_1)\big)$.
By inductive hypothesis there are $\delta_X,\delta_Z \in \m^{3e+(f-1)\v(n)}$ such that
\[
    (X_{f+1}:1:Z_{f+1}) = n (X_{f}:1:Z_{f}) = n(n^{f}X_1 + \delta_X:1:n^{f}Z_1 + \delta_Z).
\]
Since $\v(n) \leq 2e$, we have $f\v(n)+e \leq 3e+(f-1)\v(n)$, which implies 
\[
    n^{f}X_1 + \delta_X \in \m^{e+f\v(n)}, \quad n^{f}Z_1 + \delta_Z \in \m^{e+f\v(n)}.
\]
Thus, by applying $n$ times part \ref{I}, we get
\[
    X_{f+1} \equiv n^{f+1}X_1 + n\delta_X \bmod \m^{3\big(e+f\v(n)\big)}, \quad 
    Z_{f+1} \equiv n^{f+1}Z_1 + n\delta_Z \bmod \m^{3\big(e+f\v(n)\big)}.
\]
We notice that both $n\delta_X$ and $n\delta_Y$ belong to $\m^{3e+f\v(n)}$, and that $3e+f\v(n) \leq 3\big(e+f\v(n)\big)$, hence
\[
    X_{f+1} \equiv n^{f+1}X_1 \bmod \m^{3e+f\v(n)}, \quad Z_{f+1} \equiv n^{f+1}Z_1 \bmod \m^{3e+f\v(n)}.
\]

\ref{IV} Let $\m = (p)$ with $p \in \Z$. Since $\m$ is maximal, then $p$ is prime and write $\a = u p^{f}$ with $\gcd(u,p)=1$ and $f = \v(\a)$.
By the power-associativity of $\L$ (Proposition \ref{prop:PowAss}) we have
\[
    (X_{\a}:1:Z_{\a}) = u \big( p^{f} (X_1:1:Z_1) \big).
\]
As $\v(p) = 1 \leq 2e$, by part \ref{III} we have
\[
    p^{f} (X_1:1:Z_1) = ( X_{f} : 1 : Z_{f} ), \text{ with } \begin{cases}
        X_{f} \equiv p^fX_1 \bmod \m^{3e+f-1}, \\
        Z_{f} \equiv p^fZ_1 \bmod \m^{3e+f-1}.
    \end{cases}
\]
Since $\v(u) = 0$ and $X_{f}, Z_{f} \in \m^{f+e}$, by applying again part \ref{III} we have
\[
    u ( X_{f} : 1 : Z_{f} ) = (X_{\a}:1:Z_{\a}), \text{ with } \begin{cases}
        X_{\a} \equiv u X_{f} \bmod \m^{3(e+f)}, \\
        Z_{\a} \equiv u Z_{f} \bmod \m^{3(e+f)}.
    \end{cases}
\]
The conclusion follows by observing that $3e+f-1 \leq 3(e+f)$.
\endproof

Proposition \ref{prop:tech} is rather technical, and will be employed in the proofs of the next sections.
In particular, its parts \ref{III} and \ref{IV} give sufficient conditions to establish when the characteristic of the base ring equals the order or points at infinity, which does not hold in general (e.g. see \cite[Example 11.1]{SalaTaufer}).

The following results show that, differently from the affine case (Proposition~\ref{prop:AffineIntersection}), the layers of an elliptic loop $\L$ do not cover the whole $\L^{\infty}$, and they may have a non-trivial intersection at infinity, e.g. when $R$ has zero-divisors.

\begin{lemma} \label{lem:infinityValuation}
    Let $R$ be a Hausdorff discrete valued ring, and let $L_t$ be a layer of an elliptic loop $\L_{A,B}(R)$.
    Then for every $(X:1:Z) \in L_t^{\infty}$, either $Z = X = 0$ or
    \[
        \v(Z) > \v(X).
    \]
\end{lemma}
\proof Since $(X:1:Z)$ satisfies the layer equation, we have
\[
    \v(Z) \geq \min\{ \v(X^3),\v(AXZ^2),\v(BZ^3),\v(3AtX^2Z), \v(3tX),\v(9BtXZ^2),\v(A^2tZ^3)\}.
\]
Since $X,Z,t \in \m$, they have all a positive valuation.
If $Z=0$, since $R$ is Hausdorff either $X = 0$ or $\v(Z) > \v(X)$.
If $Z \neq 0$, the above inequality gives
\[
    \v(Z) \geq \min\{ 3\v(X),\v(t)+\v(3X)\},
\]
which implies $\v(Z) > \v(X)$.
\endproof

\begin{prop} \label{prop:InfIntersection}
Let $R$ be a Hausdorff discrete valued ring, and let $\L_{A,B}(R)$ be an elliptic loop. Then
\[
    \bigcap_{t \in \m} L_t^{\infty} = \begin{cases}
        \{ \O \} & \text{if } \Nil \in \{1,\infty\}, \\
        \{ (X:1:0) \}_{X \in \m^{\Nil-1}} & \text{otherwise}.
    \end{cases} 
\]
\end{prop}
\proof
The case $\Nil = 1$ is trivial, since when $R$ is a field there is only one layer $L_0 = E_{A,B}(R)$, for which $\O$ is the unique point at infinity.

When $R$ is not a field, the chain $\{\m^i\}_{i \in \Z_{\geq 0}}$ cannot stabilize at $\m \neq \< 0 \>$, hence there is an element $t \in \m \setminus \m^2$.
An infinity point $P = (X:1:Z)$ that belongs to both $L_t$ and $L_{2t}$ needs to satisfy both their layer equations, namely
\[ \begin{cases}
F(P) = t \H_F(P),\\
F(P) = 2t \H_F(P),
\end{cases} \implies
\begin{cases}
F(P) = 0,\\
t \H_F(P) = 0,
\end{cases}\implies
\begin{cases}
Z = X^3+AXZ^2+BZ^3,\\
3t X = t AZ^3 - 3t AX^2Z - 9t BXZ^2.
\end{cases}\]
By Lemma \ref{lem:infinityValuation} we have $\v(Z) \geq \v(X)$, which by means of the second equation implies
\begin{align*}
    \v(t X) &= \v(t AZ^3 - 3t AX^2Z - 9t BXZ^2) \geq \min\{ \v(t AZ^3), \v(3t AX^2Z), \v(9t BXZ^2), \Nil \} \\
    &\geq \min\{ \v(t X) + \v(AZ^2), \v(t X) + \v(AXZ), \v(t X) + \v(BZ^2), \Nil \}.
\end{align*}
This is only possible if $1+\v(X) = \v(t X) = \infty$, i.e. $t X = 0$.

If $\Nil = \infty$, then $\v(X) = \infty$, so the above relations imply $X = Z = 0$, hence $P = \O$.

If $\Nil$ is finite, then $1+\v(X) = \infty$ implies $\v(X) \geq \Nil - 1$, i.e. $X \in \m^{\Nil-1}$.
By Lemma \ref{lem:infinityValuation} we also have $\v(Z) = \Nil$, which implies $Z = 0$.
On the other side, it is easy to see that if $X \in \m^{\Nil-1}$, then for every $t \in \m$ both $F(X:1:0) = X^3$ and $t \H_F(X:1:0) = 3tX$ belong to $\m^{\Nil} = \<0\>$, hence $(X:1:0)$ belongs to every layer of $\L_{A,B}(R)$.
\endproof

\begin{remark}
By Proposition \ref{prop:tech}-\ref{I}, the set $\{ (X:1:0) \}_{X \in \m^{\Nil - 1}}$ endowed with the loop addition law is actually a group.
Thus, in absence of $3$-torsion points and with a finite nilpotency of $\m$, the intersection of all layers consists of a group at infinity, which is isomorphic to $(\m^{\Nil-1},+)$.
\end{remark}

 
 \section{More on associativity} \label{sec:MoreAss}

 In this section, we present special results about associativity with additional hypotheses on the base ring.
 First, we show that complete associativity may almost never occur among all the points of an elliptic loop.
 Afterward, we show that weak forms of associativity hold if the underlying ring has small nilpotency.
 
 
 \subsection{Non-group guarantee} \label{sec:MoreAssNO}
 
 Here we certify that elliptic loops are never groups if $\Nil \geq 3$, or if $\Nil \geq 2$ and they have at least one point of order coprime to $3$.
 We require a technical condition on the base ring, namely it needs to have at least one element whose power ideals do not immediately stabilize. This is a very broad assumption for general rings, but we remark that it is never the case if $R$ is a field.
 This is not surprising, as elliptic loops over fields agree with their underlying elliptic curves, which are associative.
 
 \begin{lemma} \label{lem:notGroupA}
     Let $R$ be a ring with an element $p \in R$ such that $p^2 \not\in \< p^3 \>$, let $\L = \L_{A,B}(R)$ be an elliptic loop and $P \in \L^a$. 
     Then
     \[
        \big( P + (p:1:p) \big) + (0:1:p) \neq P + \big( (p:1:p) + (0:1:p) \big).
     \]
 \end{lemma}
 \proof Let $P = (X:Y:1)$, and compute
 \[
    (S_1:S_2:S_3) = \big( P + (p:1:p) \big) + (0:1:p), \quad (T_1:T_2:T_3) = P + \big( (p:1:p) + (0:1:p) \big).
 \]
 Let us assume by contradiction that they are equal, then $S_iT_j-S_jT_i = 0$ for every $1 \leq i \neq j \leq 3$, hence the ideal
 \[
    I = \< p^3, S_1T_2-S_2T_1, S_1T_3-S_3T_1, S_2T_3-S_3T_2 \>
 \]
 is equal to $\< p^3 \>$. However, we can symbolically verify that
 \[ 2(4A^3+27B^2)^2p^2Y^3 \in I. \]
 Since the three quantities $2$, $\Delta_{A,B}=4A^3+27B^2$ and $Y$ are invertible in $R$ by definition of $\L$, it follows that $p^2 \in \< p^3 \>$, contradicting the hypothesis.
 \endproof
 
 \begin{lemma} \label{lem:notGroupB}
     Let $R$ be a ring with an element $p \in R$ such that $p \not\in \< p^2 \>$.
     Let also $\L = \L_{A,B}(R)$ be an elliptic loop and $P = (X:Y:1) \in \L^a$ such that $3\pi(P) \ne \O$. 
     Then
     \[
        \big( P + (X:Y+p:1) \big) + (0:1:p) \neq P + \big( (X:Y+p:1) + (0:1:p) \big).
     \]
 \end{lemma}
 \proof We straightforwardly compute
 \[
    (S_1:S_2:S_3) = \big( P + (X:Y+p:1) \big) + (0:1:p), \quad (T_1:T_2:T_3) = P + \big( (X:Y+p:1) + (0:1:p) \big).
 \]
 By inspecting the $2$-minors $S_iT_j-S_jT_i$ modulo $p^2$, we notice that they remarkably have many common factors. In fact, by denoting
 \begin{align*}
    F_1 &= A^2-3AX^2-9BX-3XY^2, \\
    F_2 &= A^3+6A^2X^2+6ABX-3AX^4+9B^2-18BX^3-Y^4, \\
    G_1 &= 10A^4X+9A^3B+2A^3Y^2-30A^2BX^2+6A^2X^5+6A^2X^2Y^2+45AB^2X+45ABX^4 \\
       & \quad +9ABXY^2+54B^3+135B^2X^3+18B^2Y^2-9BX^3Y^2, \\
    G_2 &= 2A^4-15A^2BX+30A^2X^4+6A^2XY^2+9AB^2+90ABX^3+3ABY^2-6AX^3Y^2+135B^2X^2 \\
       & \quad -27BX^5-27BX^2Y^2,
 \end{align*}
 we computationally verify that
 \[
    S_1T_2 - S_2T_1 \equiv 2pY^2F_1F_2G_1 \bmod \< p^2 \>, \quad S_2T_3 - S_3T_2 \equiv 2pY^2F_1F_2G_2 \bmod \< p^2 \>.
 \]
 Therefore, if $(S_1:S_2:S_3) = (T_1:T_2:T_3)$, one of the following should hold:
 \[
    (\text{I})\ F_1 \in \m, \quad (\text{II})\ F_2 \in \m, \quad (\text{III})\ G_1,G_2 \in \m.
 \]
 Let $F$ be the defining polynomial of $\L$. We prove that none of the above cases may occur.
 
 [Case I] We compute
 \[
    (X_3:Y_3:Z_3) = 3P
 \]
 and we verify that
 \[
    X_3, Z_3 \in \< F(P),F_1 \>.
 \]
 Since $F(P) \in \m$ by definition of $\L$, if $F_1 \in \m$ then $\pi(3P) = \O$, contradicting the hypothesis.
 
 [Case II] We compute
 \[
    (X_4:Y_4:Z_4) = 4P
 \]
 and we verify that
 \[
    X_4, Z_4 \in \< F(P),F_2 \>.
 \]
 Since $F(P) \in \m$, and $\L$ does not have points of even order by definition, then $F_2 \not\in \m$.
 
 [Case III] We straightforwardly verify that
 \[
    864Y^{10}(X^3-Y^2),\ 288Y^8(B-2X^3+2Y^2) \in \< F(P), G_1, G_2 \>,
 \]
 therefore if both $G_1,G_2 \in \m$, then we would have
 \[
    B,\, AX \in \< X^3-Y^2, B-2X^3+2Y^2, X^3+AX+B-Y^2 \> \se \m.
 \]
 This may not happen, since both $X$ and $\Delta_{A,B}$ are units of $R$.
 \endproof
 
 The points at infinity tend to associate more, as their $X$ and $Z$ entries are non-units.
 However, the next lemma shows that they do not associate as soon as $\Nil \geq 6$.
 
 \begin{lemma} \label{lem:infNONass}
    Let $R$ be a ring with an element $p \in R$ such that $p \not\in \< p^5 \>$, and let $\L = \L_{A,B}(R)$ be an elliptic loop.
    Then 
    \[
        \big( (p:1:0) + (0:1:p) \big) + (0:1:p) \neq (p:1:0) + \big( (0:1:p) + (0:1:p) \big).
    \]
 \end{lemma}
 \proof  We compute
 \[
    (S_1:S_2:S_3) = \big( (p:1:0) + (0:1:p) \big) + (0:1:p), \quad (T_1:T_2:T_3) = (p:1:0) + \big( (0:1:p) + (0:1:p) \big),
 \]
 and we define
 \[
    I = \< p^6, S_1T_2-S_2T_1, S_1T_3-S_3T_1, S_2T_3-S_3T_2 \>.
 \]
 We notice that
 \[
    972 p^5B^3(B-2), \
    36 p^5B(4A-9B^2+24B), \
    6 p^5A(2A+3B) \in I.
 \]
 If we had $(S_1:S_2:S_3) = (T_1:T_2:T_3)$, then $I = \< p^6 \>$.
 This would imply
 \[
    B(B-2), \ B(4A-9B^2+24B), \ A(2A+3B) \in \m,
 \]
 which is not possible since $\Delta_{A,B} \in R^*$ by definition of $\L$.
 \endproof
 
 
 \subsection{Low nilpotency} \label{sec:lownil}

 We now examine associativity properties that hold over Hausdorff rings with low values of $\Nil$.
 
 \begin{lemma} \label{lem:infASS}
    Let $R$ be a Hausdorff ring with $\Nil \leq 5$ and $\L$ be an elliptic loop.
    Then $\L^{\infty}$ is an abelian group.
 \end{lemma}
 \proof Given three points $\{(X_i:1:Z_i)\}_{i \in \{1,2,3\}} \se \L^{\infty}$, we explicitly compute
 \begin{align*}
    (S_1:S_2:S_3) &= \big( (X_1:1:Z_1) + (X_2:1:Z_2) \big) + (X_3:1:Z_3), \\
    (T_1:T_2:T_3) &= (X_1:1:Z_1) + \big( (X_2:1:Z_2) + (X_3:1:Z_3) \big),
 \end{align*}
 and we verify that 
 \[
    I_2 \left( \begin{bmatrix} S_1 & S_2 & S_3 \\
    T_1 & T_2 & T_3
    \end{bmatrix}\right) \se \m^5,
 \]
 which concludes the proof since $\Nil \leq 5$ implies $\m^5 = \< 0 \>$.
 \endproof
 
 By virtue of Lemma \ref{lem:infNONass}, the associativity provided by Lemma \ref{lem:infASS} is the best result we can achieve on $\L^{\infty}$ in terms of $\Nil$.
 For even smaller values of $\Nil$, this group is isomorphic to two copies of the additive group of $\m$.
 
 \begin{lemma} \label{lem:TrivialInf}
    Let $R$ be a Hausdorff ring with $\Nil \leq 3$ and $\L$ be an elliptic loop. Then
    \[
        \L^{\infty} \to (\m,+)^2, \quad (X:1:Z) \to (X,Z),
    \]
    is a well-defined group isomorphism.
 \end{lemma}
 \proof Since $\m^3 = \< 0 \>$, it immediately follows by Proposition \ref{prop:tech}-\ref{I}.
 \endproof
 
 We now provide a family of weak associativity results, which only hold for $\Nil \leq 2$.
 
 \begin{lemma}
    Let $R$ be a Hausdorff ring with $\Nil \leq 2$ and $\L$ be an elliptic loop.
    For every $P \in \L$ and $Q, R \in \L^{\infty}$, we have
    \[
        P+(Q+R) = (P+Q)+R.
    \]
 \end{lemma}
 \proof We formally compute
 \[
    (S_1:S_2:S_3) = P+(Q+R), \quad (T_1:T_2:T_3) = (P+Q)+R,
 \]
 and we verify that all the terms of
 \[
    I_2 \left( \begin{bmatrix} S_1 & S_2 & S_3 \\
    T_1 & T_2 & T_3
    \end{bmatrix}\right)
 \]
 are divisible by dregree-$2$ terms in the $X,Z$-entries of $Q$ and $R$, hence they belong to $\m^2 = \< 0 \>$.
 \endproof
 
 \begin{lemma} \label{lem:AssWithInf}
    Let $R$ be a Hausdorff ring with $\Nil \leq 2$ and $\L$ be an elliptic loop.
    For every $P,Q \in \L$ such that $\pi(P) = \pi(Q)$ and every $R_1, R_2 \in \L^{\infty}$, we have
    \[
        (P+R_1)-(Q+R_2) = (P-Q)+(R_1-R_2).
    \]
 \end{lemma}
 \proof By hypothesis there are $m_1,m_2,m_3,X_1,Z_1,X_2,Z_2 \in \m$ such that if $P = (X:Y:Z)$, then
 \[
    Q = (X+m_1:Y+m_2:Z+m_3), \quad R_1 = (X_1:1:Z_1), \quad R_2 = (X_1:1:Z_2).
 \]
 By computing
 \[
    (S_1:S_2:S_3) = (P+R_1)-(Q+R_2), \quad (T_1:T_2:T_3) = (P-Q)+(R_1-R_2),
 \]
 we straightforwardly verify that
 \[
    I_2 \left( \begin{bmatrix} S_1 & S_2 & S_3 \\
    T_1 & T_2 & T_3
    \end{bmatrix}\right) \se \< m_1,m_2,m_3,X_1,Z_1,X_2,Z_2 \>,
 \]
 from which the thesis follows by observing that $\< m_1,m_2,m_3,X_1,Z_1,X_2,Z_2 \>^2 \se \m^2 = \< 0 \>$.
 \endproof
 
 \begin{lemma} \label{lem:3PointsOverTheSame}
    Let $R$ be a Hausdorff ring with $\Nil \leq 2$ and $\L$ be an elliptic loop.
    For every $P,Q,R \in \L$ such that $\pi(P) = \pi(Q) = \pi(R)$, we have
    \[
        (P+Q)-R = P+(Q-R).
    \]
 \end{lemma}
 \proof By hypothesis there are $m_1,m_2,m_3,m_4,m_5,m_6 \in \m$ such that if $P = (X:Y:Z)$, then
 \[
    Q = (X+m_1:Y+m_2:Z+m_3), \quad R = (X+m_4:Y+m_5:Z+m_6). 
 \]
 We compute
 \[
    (S_1:S_2:S_3) = (P+Q)-R, \quad (T_1:T_2:T_3) = P+(Q-R),
 \]
 and we verify that
 \[
    I_2 \left( \begin{bmatrix} S_1 & S_2 & S_3 \\
    T_1 & T_2 & T_3
    \end{bmatrix}\right) \se \< m_1,m_2,m_3,m_4,m_5,m_6 \>^2,
 \]
 and since $\< m_1,m_2,m_3,X_1,Z_1,X_2,Z_2 \>^2 \se \m^2 = \< 0 \>$, then the triple $P,Q,-R$ is associative.
 \endproof
 
 \begin{remark}
    In the setting of Lemma \ref{lem:3PointsOverTheSame} we can avoid useless parentheses: as the result is independent of the association order, we simply denote it by $P+Q-R$.
 \end{remark}
 
 \begin{prop} \label{prop:HeavyProp}
    Let $R$ be a Hausdorff ring with $\Nil \leq 2$ and $\L$ be an elliptic loop.
    For every $P_1,P_2,P_3,Q_1,Q_2,Q_3 \in \L$ such that $\pi(P_1) = \pi(P_2) = \pi(P_3)$ and $\pi(Q_1) = \pi(Q_2) = \pi(Q_3)$, we have
    \[
        (P_1+P_2-P_3) + (Q_1+Q_2-Q_3) = (P_1+Q_1)+(P_2+Q_2)-(P_3+Q_3).
    \]
 \end{prop}
 \proof Since $\L$ has no even-order points, we may assume $P_1 = (X_1:1:Z_2)$, $Q_1 = (X_2:1:Z_2)$ and we can consider $m_1,m_2,m_3,m_4,m_5,m_6,m_7,m_8 \in \m$ such that
 \begin{align*}
    &P_2 = (X_1+m_1:1:Z_1+m_2), &&P_3 = (X_1+m_3:1:Z_1+m_4), \\
    &Q_2 = (X_2+m_5:1:Z_2+m_6), &&Q_3 = (X_2+m_7:1:Z_2+m_8).
 \end{align*}
 We formally compute
 \[
    (S_1:S_2:S_3) = (P_1+P_2-P_3) + (Q_1+Q_2-Q_3), \quad
    (T_1:T_2:T_3) = (P_1+Q_1)+(P_2+Q_2)-(P_3+Q_3).
 \]
 An intensive computation shows that
 \[
    I_2 \left( \begin{bmatrix} S_1 & S_2 & S_3 \\
    T_1 & T_2 & T_3
    \end{bmatrix}\right) \se \< m_1,m_2,m_3,m_4,m_5,m_6,m_7,m_8 \>^2,
 \]
 which concludes the proof since $\< m_1,m_2,m_3,m_4,m_5,m_6,m_7,m_8 \>^2 \se \m^2 = \< 0 \>$.
 \endproof
 
 As a consequence of Proposition \ref{prop:HeavyProp}, we observe that scalar multiplication is distributive on the associative triples of Lemma \ref{lem:3PointsOverTheSame}.
 
 \begin{lemma} \label{lem:AssMultiples}
    Let $R$ be a Hausdorff ring with $\Nil \leq 2$ and $\L$ be an elliptic loop.
    For every triple $P_1,P_2,P_3 \in \L$ such that $\pi(P_1) = \pi(P_2) = \pi(P_3)$ and every $m \in \Z$, we have
    \[
        m(P_1 + P_2 - P_3) = mP_1 + mP_2 - mP_3.
    \]
 \end{lemma}
 \proof We prove it by induction on $m$.
 
 [$m = 1$] There is nothing to prove. 
 
 [$m \to m+1$] By inductive hypothesis we have
 \[
    (m+1)(P_1 + P_2 - P_3) = (P_1 + P_2 - P_3) + (mP_1 + mP_2 - mP_3).
 \]
 Since $\pi(mP_1) = \pi(mP_2) = \pi(mP_2) = m\pi(P_1)$, then the thesis follows by Proposition \ref{prop:HeavyProp}.
 \endproof
 
 \section{Finite order in the same fiber} \label{sec:PointsFiniteOrd}
 
 In this section, we employ the results of Section \ref{sec:lownil} for characterizing the torsion elements of elliptic loops when $\Nil \leq 2$.
 
 For every integer $q \in \Z$ and elliptic loop $\L$, we denote the points of $\L$ annihilated by the multiplication-by-$q$ map as 
 \[
    \L_q = \{ P \in \L \ | \ qP = \O \}.
 \]
 In other terms, the set $\L_q$ is made of the points of $\L$ whose order is finite and divides $q$.
 Moreover, given a point $P \in \L_q$, we denote the points of $\L_q$ lying over $\pi(P)$ by
 \[
    \L_{q/P} = \{ Q \in \L_q \ | \ \pi(Q) = \pi(P) \},
 \]
 and we denote the set of their differences by
 \[
    \D_{q/P} = \{P_1 - P_2\}_{P_1, P_2 \in \L_{q/P}}.
 \]
 
 \begin{thm} \label{thm:qTors}
    Let $R$ be a Hausdorff ring with $\Nil \leq 2$ and $\L$ be an elliptic loop.
    For every $q \in \Z$ and $P \in \L$ we have that $\D_{q/P}$ is a subgroup of $\L^{\infty}$, and
    \[
        \L_{q/P} = P + \D_{q/P}. 
    \]
 \end{thm}
 \proof By Lemma \ref{lem:infASS} we know that $\L^{\infty}$ is an abelian group, which clearly contains $\D_{q/P}$.
 Moreover, by Lemma \ref{lem:AssWithInf}, for every $P_1,P_2,P_3,P_4 \in \L_{q/P}$ we have
 \[
    (P_1 - P_2) + (P_3 - P_4) = (P_1+P_3-P_4) - P_2,
 \]
 and $P_1+P_3-P_4 \in \L_{q/P}$ by Lemma \ref{lem:AssMultiples}.
 Hence the addition law is closed on $\D_{q/P}$, and since $\O = P_1 - P_1 \in \D_{q/P}$ and $-(P_1-P_2) = P_2 - P_1 \in \D_{q/P}$, then $\D_{q/P}$ is a subgroup of $\L^{\infty}$.
 
 The inclusion $\L_{q/P} \se P + \D_{q/P}$ is clear, since every $Q \in \L_{q/P}$ may be written as $Q = P + (Q - P)$ by Lemma \ref{lem:WeakAssociativity}.
 On the other side, for every $P_1, P_2 \in \L_{q/P}$, by Lemma \ref{lem:AssMultiples} we have
 \[
    q( P + P_1 - P_2 ) = qP + qP_1 - qP_2 = 0,
 \]
 and since $\pi$ is linear, also $\pi(P + P_1 - P_2) = \pi(P)$.
 Thus, we conclude $P + P_1 - P_2 \in \L_{q/P}$.
 \endproof
 
 \begin{prop} \label{prop:qLine}
 Let $R$ be a Hausdorff ring with $\Nil \leq 2$ and $\L = \L_{A,B}(R)$ be an elliptic loop.
 For every $q \in \Z$ and $P \in \L$, if $\D_{q/P}$ is cyclic then there is a projective line $L \se \P^2(R)$ such that
 \[
    \L_{q/P} \se L.
 \]
 Moreover, if also $\m$ is a principal $\Z$-module, then there is a projective line $L' \se \P^2(R)$ such that
 \[
    \L_{q/P} = L' \cap \pi^{-1}\big(\pi(P)\big).
 \]
 \end{prop}
 \proof Let $(X:1:Z) \in \L_{q/P}$ and $(m_x:1:m_z) \in \L^{\infty}$ be a generator of $\D_{q/P}$.
 By Lemma \ref{lem:TrivialInf} and Theorem \ref{thm:qTors}, there is $(m_x:1:m_z) \in \D_{q/P}$ such that
 \[
    \L_{q/P} = \{ (X:1:Z) + (km_x:1:km_z) \}_{k \in \Z}.
 \]
 By defining
 \begin{align*}
    \a &= A^2Z^2m_z - AX^2m_z - 2AXZm_x - 6BXZm_z - 3BZ^2m_x + m_x, \\
    \b &= 2AXZm_z + AZ^2m_x + 3BZ^2m_z + 3X^2m_x + m_z.
 \end{align*}
 we can straightforwardly verify that, since $\m^2 = \< 0 \>$, then
 \[
    (X:1:Z) + (km_x:1:km_z) = (X+k\a:1:Z+k\b).
 \]
 Hence, all the points of $\L_{q/P}$ belong to the projective line $L$ defined by
 \[
    -\b x + (\b X - \a Z) y + \a z = 0.
 \]
 As for the moreover part, we notice that both $\a,\b \in \m$, hence if $\m = \g \Z$, then we can write $\a = \a' \g,\ \b = \b' \g$ and check that again we have
 \[
    -\b' x + (\b' X - \a' Z) y + \a' z = 0,
 \]
 so $\L_{q/P}$ is contained in the line $L'$ defined by the above equation (and in $\pi^{-1}\big(\pi(P)\big)$ by definition).
 On the other side, for every $Q \in \pi^{-1}\big(\pi(P)\big)$ there are $s_x,s_z \in \Z$ such that $Q = (X+s_x \g : 1 : Z+s_z \g)$.
 If also $Q \in L'$, then
 \[
    0 = -\b' (X+s_x \g) + (\b' X - \a' Z) + \a' (Z+s_z \g) = -\b's_x \g + \a' s_z \g.
 \]
 If $\a',\b' \in \m$, then $Q = P \in \L_{q/P}$.
 Otherwise, one of them is invertible, without losing of generality $\a' \in R^*$.
 In this case we have
 \[
    X+s_x \g = X + \frac{s_x}{\a'} \a' \g = X + \frac{s_x}{\a'} \a, \quad Z+s_z \g = Z + \frac{s_x}{\a'} \b' \g = Z + \frac{s_x}{\a'} \b.
 \]
 Hence, by defining $k = \frac{s_x}{\a'} \in \Z$, we conclude
 \[
    Q = (X + k\a : 1 : Z + k\b) = (X:1:Z) + k(m_x:1:m_z) \in \L_{q/P},
 \]
 which proves the second inclusion.
 \endproof
 
 
\section{The $\Z/p^e\Z$-case} \label{sec:ZNZ}

 In this section, we consider an integer prime $p$ and a positive integer $e$, and we examine elliptic loops over the ring $R=\Z/p^e\Z$, whose maximal ideal is $\m = \< p \>$.
 The $\m$-adic discrete valuation of this ring is usually referred to as the $p$-adic valuation, and it is denoted by $\v_p$.
 
 An elliptic loop $\L$ over this ring projects to an elliptic curve $E$ defined over $\F_p$, whose possible group structures are known \cite{Ruck,Voloch}, although they may be difficult to be effectively computed.
 
 We denote by $q = |E|$ the size of the underlying curve.
 In this setting, $\pi$ is a $p^{2(e-1)}$-covering, thus we have
 \[
    |\L^{\infty}| = p^{2(e-1)}, \quad |\L^{a}| = (q-1)p^{2(e-1)}.
 \]
 By Proposition \ref{prop:loop} and \ref{prop:PowAss} we know that such $\L$ is a power-associative abelian algebraic loop, while by Lemma \ref{lem:notGroupA} and \ref{lem:notGroupB} we also know that it is not a group if one of the following conditions hold:
 \[
    (\text{A})\ e \geq 3 \quad (\text{B})\ e \geq 2 \text{ and } p \geq 17.
 \]
 In fact, when $p \geq 17$ the base curve always has a point of order other than $3$, since it needs to have at least $10$ elements by the Hasse Theorem \cite[Theorem 4.2]{Washington} and its group of points has rank at most $2$ \cite[Theorem 4.1]{Washington}.
 
 By testing every existent elliptic curve defined over a prime field $\F_p$ (with $p \leq 17$) with a group of points isomorphic to either $\Z/3\Z$ or $\Z/3\Z \times \Z/3\Z$, one may verify \cite{TauThesis} that there are precisely $6$ elliptic loops over local quotients of $\Z$ that are actual groups, namely
 \begin{align*}
    \L_{4,2}(\Z/25\Z) \simeq \L_{4,3}(\Z/25\Z) &\simeq \Z/5\Z \times \Z/15\Z,\\
    \L_{0,4} (\Z/49\Z) &\simeq \Z/7\Z \times \Z/21\Z,\\
    \L_{0,2} (\Z/49\Z) &\simeq \Z/21\Z \times \Z/21\Z,\\
    \L_{0,3}(\Z/169\Z) \simeq \L_{0,10}(\Z/169\Z) &\simeq \Z/39\Z \times \Z/39\Z.
 \end{align*}
 We also know that $\L^{a}$ is (set-theoretically) made of all the affine points of elliptic curves lifting $E(\F_p)$ (Lemma \ref{lem:AffinePointsofLoop}), and when $3 \nmid q$ it is stratified by the layers of $\L$ (Proposition \ref{prop:AffineIntersection}).
 
 Over these special rings, the results of Section \ref{sec:InfinityPart} lead to a complete characterization of the infinity part $\L^{\infty}$.
 
 \begin{thm} \label{thm:InfStructure}
    Let $\L = \L_{A,B}(\Z/p^e\Z)$ be an elliptic loop. Then
    \begin{enumerate}[label=$(\roman*)$]
        \item \label{i1} $\langle (p:1:0) \rangle \cap \langle (0:1:p) \rangle = \{ \O \}$.
        \item \label{i2} $| \langle (p:1:0) \rangle | = | \langle (0:1:p) \rangle | = p^{e-1}$.
        \item \label{i3} For every $P \in \L^{\infty}$ there are uniquely determined integers $0 \leq \a,\b \leq p^{e-1}-1$, such that
        \[
            P = \a (p:1:0) + \b (0:1:p).
        \]
    \end{enumerate}
 \end{thm}
 \proof 
 \ref{i1} By Proposition \ref{prop:tech}-\ref{IV}, every $(X:1:Z) = \a (p:1:0)$ satisfies
 \[
    Z \equiv 0 \mod p^{3\v_p(X)-1} \implies X = Z = 0 \text{ or } \v_p(Z) > \v_p(X). 
 \]
 Similarly, if the same point may also be written as $\b (0:1:p)$, then
 \[
    X \equiv 0 \mod p^{3\v_p(Z)-1} \implies X = Z = 0 \text{ or } \v_p(X) > \v_p(Z).
 \]
 The unique point satisfying both the above conditions is $(X:1:Z) = \O$.
 
 \ref{i2} By Proposition \ref{prop:tech}-\ref{IV} we know that if $(X:1:Z) = \a (p:1:0)$, then
 \[
    X \equiv \a p \bmod p^{\v_p(\a) + 2},
 \]
 hence the lowest possible $\a \in \Z_{\geq 1}$ such that $X = 0$ is $\a = p^{e-1}$, which is therefore the order of the group $\langle (p:1:0) \rangle$.
 In the same way, we prove the analogous result on $\langle (0:1:p) \rangle$.
 
 \ref{i3} The case $e=1$ is trivial since $\L_{A,B}(\F_p) = E_{A,B}(\F_p)$ has a unique point at infinity, hence we may assume $e \geq 2$.
 Since $|\L^{\infty}| = |\ker \pi| = p^{2(e-1)}$, it is sufficient to prove that for every $\a_1,\a_2,\b_1,\b_2 \in \Z$ and every $2 \leq \e \leq e$ we have
 \[
    \a_1(p:1:0) + \b_1(0:1:p) = \a_2(p:1:0) + \b_2(0:1:p) \implies \begin{cases}
    \a_1 \equiv \a_2 \bmod p^{\e-1},\\
    \b_1 \equiv \b_2 \bmod p^{\e-1}.
    \end{cases}
 \]
 We prove it by induction on $\e \in \Z_{\geq 2}$.
 
 [$\e=2$] It follows from Proposition \ref{prop:tech}-\ref{I}.
 
 [$\e \to \e+1$] By inductive hypotheses there are $\a,\b \in \Z$ such that we have $\a_2 = \a_1 + \a p^{\e-1}$ and $\b_2 = \b_1 + \b p^{\e-1}$.
 Since $\L$ is power-associative, we can write
 \[
    \a_2(p:1:0) + \b_2(0:1:p) = \big( \a_1(p:1:0)+\a p^{\e-1}(p:1:0) \big) + \big( \b_1(0:1:p)+\b p^{\e-1}(0:1:p) \big).
 \]
 Moreover, if we define
 \[
    (\delta_{1}:1:\delta_{2}) = \a p^{\e-1}(p:1:0), \quad (\delta_3:1:\delta_4) = \b p^{\e-1}(0:1:p),
 \]
 by Proposition \ref{prop:tech}-\ref{IV} we have
 \[
    \delta_1 \equiv \a p^\e \bmod p^{\e+1}, \quad \delta_2 \equiv \delta_3 \equiv 0 \bmod p^{\e+1}, \quad \delta_4 \equiv \b p^\e \bmod p^{\e+1}.
 \]
 Furthermore, by defining
 \[
    (X_{\a}:1:Z_{\a}) = \a_1 (p:1:0), \ (X_{\b}:1:Z_{\b}) = \b_1 (0:1:p), \ (X:1:Z) = (X_{\a}:1:Z_{\a}) + (X_{\b}:1:Z_{\b})
 \]
 and
 \[
    (X_{1}:1:Z_{1}) = (X_{\a}:1:Z_{\a}) + (\delta_{1}:1:\delta_{2}), \quad
    (X_{2}:1:Z_{2}) = (X_{\b}:1:Z_{\b}) + (\delta_{3}:1:\delta_{4}),
 \]
 by Proposition \ref{prop:tech}-\ref{II} we have
 \[
    X_{1} \equiv X_{\a} + \a p^\e, \quad Z_{1} \equiv Z_{\a}, \quad 
    X_{2} \equiv X_{\a}, \quad Z_{2} \equiv Z_{\a} + \b p^\e.
 \]
 Thus, $\a_2(p:1:0) + \b_2(0:1:p)$ is equal to
 \[
    (X_f:1:Z_f) = (X_1:1:Z_1) + (X_2:1:Z_2),
 \]
 where by Proposition \ref{prop:tech}-\ref{II} we have
 \[
    X_f \equiv X + \a p^\e, \quad  Z_f \equiv Z + \b p^\e.
 \]
 In conclusion, if we have $(X:1:Z) = (X_f:1:Z_f)$, then $p|\a$ and $p|\b$, which implies $\a_1 \equiv \a_2 \bmod p^{\e}$ and $\b_1 \equiv \b_2 \bmod p^{\e}$.
 \endproof
    
 Theorem \ref{thm:InfStructure} shows that $\L^{\infty}$ may be set theoretically thought of as
 \[
    \L^{\infty} \longleftrightarrow \langle (p:1:0) \rangle \times \langle (0:1:p) \rangle.
 \]
 However, we remark that the above correspondence does not preserve the addition law, as the sum in $\L^{\infty}$ is not necessarily the product-sum of these subloops, unless $e \leq 5$ (Lemma \ref{lem:infASS}). 

 The next lemma shows that $\langle (0:1:p) \rangle$ has a trivial intersection with every layer.
 
 \begin{lemma} \label{lem:ForbiddenInf}
    Let $\L = \L_{A,B}(\Z/p^e\Z)$ be an elliptic loop. For every $t \in p\Z/p^e\Z$, we have
    \[
        L_t \cap \langle (0:1:p) \rangle = \{ \O \}.
    \]
 \end{lemma}
 \proof By Proposition \ref{prop:tech} every non-zero multiple $(X:1:Z) = \a (0:1:p)$ satisfies
 \[
    X \equiv 0 \bmod p^{\v_p(\a)+2}, \quad Z \equiv \a p \bmod p^{\v_p(\a)+2}.
 \]
 Hence, we have $\v_p(X) \geq \v_p(Z)$.
 If also $(X:1:Z) \in L_t$, then by Lemma \ref{lem:infinityValuation} we conclude $(X:1:Z) = \O$.
 \endproof

 Not only the infinity part of the loop may be characterized, but also the infinity part of its layers can be established.
 Indeed, all the $p^{e-1}$ layers $L_t$ of $\L$ have a cyclic infinity part.
 
 \begin{prop} \label{prop:LayStructure}
    Let $\L = \L_{A,B}(\Z/p^e\Z)$ be an elliptic loop. For every $t \in p\Z/p^e\Z$ there is a unique $Z_t \in p\Z/p^e\Z$ such that
    \[
        L_t^{\infty} = \langle (p:1:Z_t) \rangle.
    \]
 \end{prop}
 \proof We consider the polynomial $F-t\H_F(p,1,z) \in \Z[z]$.
 Since $F-t\H_F(p,1,z)(p,1,0) \equiv 0 \bmod p$ and its derivative in $z = 0$ is constantly $-1$, by the Hensel Lemma there is a unique $Z_t \in p\Z/p^e\Z$ such that $(p:1:Z_t) \in L_t$.
 Since $\pi|_{L_t}$ is a $p^{e-1}$-covering of the base curve, it is sufficient to show that the order of $(p:1:Z_t)$ is $p^{e-1}$.
 For a given $\a \in \Z$, if we define
 \[
    (X_{\a}:1:Z_{\a}) = \a (p:1:Z_t),
 \]
 then by proposition \ref{prop:tech}-\ref{IV} we have
 \[
    X_{\a} \equiv \a p \bmod p^{\v_p(\a)+2}, \quad Z_{\a} \equiv \a Z_t \bmod p^{\v_p(\a)+2}.
 \]
 Thus, the minimal positive $\a$ that makes $(X_{\a}:1:Z_{\a}) = \O$ is $\a = p^{e-1}$.
 \endproof
 
 When the size $q$ of the underling curve is coprime to $p$, the short exact sequence \eqref{eq:ses} splits, thus Proposition \ref{prop:LayStructure} provides the group structure of layers:
 \[
    L_t \simeq \Z/p^{e-1}\Z \times E_{A,B}(\F_p).
 \]
 We also remark that the entry $Z_t$ in the generator of $L_t^{\infty}$ may be explicitly computed by truncating to the correct exponent of the classical series \cite[Chapter IV]{Silverman}.
 This result extends the group classification of elliptic curves over $\Z/N\Z$ \cite[Theorem 17]{SalTau} to every layer of $\L$.
 
 Finally, we observe that over such rings, the group $\D_{q/P}$ of Section \ref{sec:PointsFiniteOrd} is almost always cyclic, hence Proposition \ref{prop:qLine} applies: there is a projective line that locates the points lying over the same base point that have $q$-torsion.
 
 \begin{prop}
 Let $\L = \L_{A,B}(\Z/p^2\Z)$ be an elliptic loop, and let $P \in \L$ be a point of order $q$, such that $\gcd(3p,q) = 1$. Then
 \[\begin{cases}
    \D_{q/P} = \{\O\}         & \text{if } P \in \L^{\infty}, \\
    \D_{q/P} \simeq \Z/p\Z    & \text{if } P \in \L^{a}.
 \end{cases}\]
 \end{prop}
 \proof By Theorem \ref{thm:qTors} we know that $\D_{q/P}$ is a subgroup of $\L^{\infty}$, hence by Theorem \ref{thm:InfStructure} it is a subgroup of $\Z/p\Z \times \Z/p\Z$.
 
 If $P \in \L^{\infty}$, then its order is divisible by both $p$ and $q$, hence $P = \O$.
 
 Conversely, if $P \in \L^{a}$ then $\pi(P) \in E_{A,B}(\F_p)$ is an affine point of order coprime to $3$.
 Since $\gcd(p,q) = 1$, by Proposition \ref{prop:LayStructure} all the layers $L_t$ of $\L$ have the same structure, namely
 \[
    L_t \simeq \Z/p\Z \times E_{A,B}(\F_p).
 \]
 Hence, $\L_{q/P}$ intersects the every layer $L_t$ in precisely one point $P_t$, i.e. the point corresponding to $\big(0,\pi(P)\big)$ under the above group isomorphism.
 Different layers never intersect over points with order coprime to $3$, as proved in Proposition \ref{prop:AffineIntersection}.
 Therefore, all the $P_t$'s are different, and the group
 \[
    \D_{q/P} = \{ P_{t} - P_{s} \}_{s,t \in p\Z/p^2\Z}
 \]
 contains at least $p$ different points, but it cannot have size $p^2$ since some elements are repeated (e.g. $P_{t}-P_{t} = P_{s}-P_{s}$).
 Thus, we must have $|\D_{q/P}| = p$.
 \endproof
 
 
 \section{Conclusions and open problems} \label{sec:conclusions}

 We presented the construction of a new family of power-associative abelian loops with a canonical projection on an elliptic curve $E$ defined over the residue field of a given local ring $(R,\m)$.
 
 When $\m = \< 0 \>$ (i.e. $R$ is a field), the results of this paper coincide with known results of elliptic curves over fields, but when $\m$ is non-trivial this construction simultaneously encodes all the lifts of the base curve.
 In particular, layers appear to correspond 1-1 to elliptic curves over $R$ with a canonical projection on $E$, but the identification of the element $t \in \m$ that gives rise to the layer $L_t$ corresponding to the canonical lift of $E$ is still open.
 
 Determining the structure of elliptic loops over rings different from $\Z/p^e\Z$ is also considered a challenging research line, as it is deeply connected to determining the structure of the formal group associated with elliptic curves over local rings.
 Besides, we identified large associative structures over $\Z/p^e\Z$, but it is still only conjectured that layers and groups of size $p^{e-1}$ at infinity are the largest possible subgroups of an elliptic loop.
 
 Finally, the investigation of the abelian varieties admitting a loop sup-structure is open, as well as the few cases excluded from the present paper ($6 \in \m$ or points with order not coprime to $6$).

 
 \section*{Acknowledgments}
 This paper has been supported by the European Union’s H2020 Programme under grant agreement number ERC-669891 and by the Research Foundation - Flanders (FWO), project 12ZZC23N.
 The authors thank Igor Semaev, Francesco Pappalardi, Michael Kinyon and the anonymous reviewers for their precious observations, and Nicolò Cangiotti for his support on the computational part.

 \section*{Statements and Declarations}
 Authors declare that there are no financial or non-financial interests that are directly or indirectly related to the work submitted for publication. \\
 Data sharing is not applicable to this article as no datasets were generated or analyzed during the current study.\\
 On behalf of all authors, the corresponding author states that there is no conflict of interest.


\begin{thebibliography}{}

\bibitem{Magma}
W. Bosma, J. Cannon, C. Playoust,
\emph{The Magma algebra system. I. The user language},
J. Symbolic Comput. 24, 1997, pp. 235--265.

\bibitem{BosLen}
W. Bosma, H. W. Lenstra,
\emph{Complete Systems of Two Addition Laws for Elliptic Curves},
J. Number Theory 53, 1995, pp. 229--240.

\bibitem{Brown}
W. C. Brown,
\emph{Matrices over commutative rings},
Marcel Dekker, 1986.

\bibitem{CilOtt}
C. Ciliberto, G. Ottaviani,
\emph{The Hessian map},
Int. Math. Res. Not., 2020, pp. 5781--5817.

\bibitem{Deuring}
M. Deuring,
\emph{Die typen der multiplikatorenringe elliptischer funktionenk\"oper},
Abh. Math. Sem. Univ. Hamburg 14, 1941, pp. 197--272.

\bibitem{Ethe}
I. M. H. Etherington,
\emph{Quasigroups and cubic curves},
Proc. Edinburgh Math. Soc. 14, 1964, pp. 273--291.

\bibitem{GaMcDO}
G. Ganske, B. R. McDonald,
\emph{Finite Local Rings},
Rocky Mt. J. Math. 3, 1973, pp. 521--540.

\bibitem{Kanevsky}
D. Kanevsky,
\emph{An example of a non-associative Moufang loop of point classes on a cubic surface},
arXiv:2104.05118, 2021.


\bibitem{ECOverSchemes}
N. M. Katz, B. Mazur,
\emph{Arithmetic Moduli of Elliptic Curves},
Ann. Math. Studies 108, Princeton University Press, 1985.

\bibitem{Lenstra}
H. W. Lenstra,
\emph{Elliptic curves and number-theoretic algorithms},
Proc. International Congress of Mathematicians, 1986, pp. 99--120.

\bibitem{Manin}
Yu. I. Manin,
\emph{Cubic Forms},
North Holland, 1986.

\bibitem{Manis}
M. E. Manis,
\emph{Valuations on a Commutative Ring},
Proc. Am. Math. Soc. 20, 1969, pp. 193--198.

\bibitem{Norman}
P. Norman,
\emph{Lifting abelian varieties},
Invent. Math. 64, 1981, pp. 431–443.

\bibitem{Ruck}
H. G. Rück,
\emph{A Note on Elliptic Curves Over Finite Fields},
Math. Comp. 49, 1987, pp. 301--304.

\bibitem{SalTau}
M. Sala, D. Taufer, 
\emph{The group structure of elliptic curves over $\Z/N\Z$},
arXiv:2010.15543, 2020.

\bibitem{SalaTaufer}
M. Sala, D. Taufer, 
\emph{A survey on the group of points arising from elliptic curves with a Weierstrass model over a ring},
Int. J. Group Theory 12, 2023, pp. 177--196.

\bibitem{SatohAraki}
T. Satoh, K. Araki,
\emph{Fermat quotients and the polynomial time discrete log algorithm for anomalous elliptic curves},
Commentarii Math. Univ. St. Pauli. 47, 1998, pp. 81--92.

\bibitem{Schafer}
R. D. Schafer,
\emph{An introduction to nonassociative algebras}, Academic Press, 1966.

\bibitem{Silverman}
J. H. Silverman, 
\emph{The Arithmetic of Elliptic Curves, 2nd Edition},
Springer, 2016.

\bibitem{Smart}
N. Smart,
\emph{The discrete logarithm on elliptic curves of trace one},
J. Cryptology 12, 1999, pp. 193--196.

\bibitem{TauThesis}
D. Taufer,
\emph{Elliptic Loops},
University of Trento, 2020. Available at \url{https://iris.unitn.it/handle/11572/265846}.

\bibitem{Voloch}
J. F. Voloch,
\emph{A note on elliptic curves over finite fields},
Bull. Soc. Math. France 116, 1988, pp. 455--458.

\bibitem{Washington}
L. C. Washington,
\emph{Elliptic curves, number theory and cryptography},
Chapman \& Hall / CRC, 2008.

\bibitem{MyCode}
The magma code for formal verification and computational testing is available on Github:
\url{https://github.com/DTaufer/Elliptic-Loops}.
\end{thebibliography}
\end{document}